\documentclass[smallcondensed]{svjour31}     % onecolumn (ditto)
\pdfoutput=1 
\smartqed  % flush right qed marks, e.g. at end of proof
\usepackage{graphicx}

\usepackage{amssymb,bm,nicefrac,mathtools} % amsmath
\usepackage{tabularx,booktabs,scalerel}
\usepackage[section]{algorithm}
\usepackage{wrapfig,float,xcolor}
\usepackage{grffile,subcaption}
\usepackage{mdframed,textpos}
\usepackage{bigints,relsize}
\usepackage{nameref}
\usepackage{todonotes}
\reversemarginpar
\usepackage{chngcntr}
\counterwithout{algorithm}{section}  %% Remove to return to section labeling on algorithms
\usepackage{tikz, tikzscale, adjustbox}
	\usetikzlibrary{shapes,arrows}
	\usetikzlibrary{external}

\usepackage{pgfplots}
	\pgfplotsset{compat=newest}
	\usetikzlibrary{plotmarks}
	\usetikzlibrary{arrows.meta}
	\usepgfplotslibrary{patchplots}
	\usepgfplotslibrary{external} 
	
	\pgfplotsset{every axis/.append style={
	label style={font=\small},
         tick label style={font=\small},
         title style={font=\bfseries\LARGE},
         grid style={line width=.1pt, draw=gray!25},
	major grid style={line width=.2pt,draw=gray!50}
	}}
	
	\pgfplotsset{every tick label/.append style={font=\tiny}}

\usepackage{hyperref,cite}
\usepackage[capitalize,noabbrev]{cleveref} 

\newcommand{\bfparp}{\mathbf{p}}
\newcommand{\cbfG}{\mathbf{\mathcal{G}}}

\newcommand{\cbfH}{\mathbf{\mathcal{H}}}
\DeclareMathOperator*{\argmin}{arg\,min}                   % arg min
\newcommand{\bfDelta}{{\boldsymbol{\Delta}}}

\newcommand{\bfpi}{{\boldsymbol{\pi}}}

\newcommand{\bfA}{{\bf A}}
\newcommand{\bfB}{{\bf B}}
\newcommand{\bfC}{{\bf C}}

\newcommand{\bfE}{{\bf E}}

\newcommand{\bfI}{{\bf I}}

\newcommand{\bfK}{{\bf K}}

\newcommand{\bfM}{{\bf M}}

\newcommand{\bfb}{{\bf b}}
\newcommand{\bfc}{{\bf c}}

\newcommand{\bfp}{{\bf p}}

\newcommand{\bfx}{{\bf x}}

\newcommand{\bfz}{{\bf z}}

\newcommand{\calA}{\mathcal{A}}

\newcommand{\calH}{\mathcal{H}}
\newcommand{\calI}{\mathcal{I}}

\newcommand{\calP}{\mathcal{P}}

\newcommand{\calR}{\mathcal{R}}

\newcommand{\bbC}{\mathbb{C}}

\newcommand{\bbV}{\mathbb{V}}

\newcommand{\qedsymbol}{\hfill$\blacksquare$}

\newcommand{\localH}{\mathfrak{h}}

\newcommand{\rmd}{\,\mathsf{d}\hspace{0.05em}}

\newcommand{\rmi}{\mathrm{i}\hspace{0.1em}}

\newif\ifuseTikzGraphs
\useTikzGraphstrue
\crefformat{equation}{(#2#1#3)}
\graphicspath{{Graphics/}}

\begin{document}

\title{Empirical least-squares fitting \\ of parametrized dynamical systems
}
\author{Alexander Grimm         \and
       Christopher Beattie \and
       Zlatko Drma\v{c} \and
        Serkan Gugercin
}

\institute{A. Grimm \at
Department of Mathematics, Virginia Tech\\
              Blacksburg, VA-24061-0123, U.S.A \\
              \email{alex588@vt.edu}           %  \\
	\and
	C. Beattie \at
Department of Mathematics, Virginia Tech\\
              Blacksburg, VA-24061-0123, U.S.A \\
            \email{beattie@vt.edu}
           \and
           S. Gugercin \at
Department of Mathematics, Virginia Tech\\
              Blacksburg, VA-24061-0123, U.S.A \\
              \email{gugercin@math.vt.edu}           %  \\
              \and
              Z. Drma\v{c} \at 
              Department of Mathematics, University of Zagreb\\
              Zagreb, Croatia\\
              \email{drmac@math.hr}
}

\maketitle

\begin{abstract}
Given a set of response observations for a parametrized dynamical system, we seek a parametrized dynamical model that will yield uniformly small response error over a range of parameter values yet has low order.  Frequently, access to internal system dynamics or equivalently, to realizations of the original system is either not possible or not practical; only response observations over a range of parameter settings might be known.  Respecting these typical operational constraints, we propose a two phase approach that first encodes the response data into a high fidelity intermediate model of modest order, followed then by a compression stage that serves to eliminate redundancy in the intermediate model.  For the first phase, we extend non-parametric least-squares fitting approaches so as to accommodate parameterized systems.
This results in a (discrete) least-squares problem formulated with respect to both frequency and parameter that identifies ``local" system response features.  The second phase uses an $\mathbf{\mathcal{H}}_2$-optimal model reduction strategy accommodating the specialized parametric structure of the intermediate model obtained in the first phase. The final compressed model inherits the parametric dependence of the intermediate model and maintains the high fidelity of the intermediate model, while generally having dramatically smaller system order.  We provide a variety of numerical examples demonstrating our approach. 
\keywords{Parametric Systems \and Data Driven Modeling \and Least Squares }
\subclass{ 35B30 \and  41A20 \and 93A15 \and 93B15 \and 93C05   \and 93E24 }
\end{abstract}

%!TEX root = ms.tex
%  +  +  +  +  +  +  +  +  +  +  +  +  +  +  +  +  +  +  +  +  +  +  +  +  +  +  +  +  +  +  +  +  +  +  +  +  +  +  +  +
\section{Introduction}\label{sec:pvf:intro}

In many areas of study in science and engineering, the dynamics that govern processes of interest are 
expected to vary with respect to a given set of parameters describing, say, viscosity, temperature distribution, reaction and flow rates, etc. 
Often these parameters may be deduced from dimensional analysis or dynamic similarity considerations while detailed dynamics of the system may be inaccessible to direct modeling.  Although access to internal dynamics may be lacking,
there may be an abundance of accurate frequency response measurements available, which may further be classified according to what system parameter values were in effect when the system responses were observed.   
The problem that we address here involves building up an \emph{empirical parsimonious parameterized dynamical system model} that fits the observed system response with high fidelity over a desired range of parameter values.  We interpret \emph{parsimonious} in this context to mean low system order.  The derived dynamical system may then be used as an efficient surrogate to predict operational behavior or to produce effective control and optimization strategies.

\subsection{Basic problem formulation}
Since linear time invariant dynamical systems have a frequency domain response representation that involves rational functions, a natural formulation for the task at hand leads one to a data fitting problem using rational functions; this will be our principal focus.  For simplicity, we assume that the original dynamics vary with respect to a single parameter $p$, and arise as a single-input/single-output (SISO) linear time-invariant system associated with a ($p$-dependent) transfer function, $\cbfH(s,p)$, that is unknown but is nonetheless accessible to sampling in the sense that for a variety of operating conditions distinguished by parameter values $p=\mu_1,\ldots, \mu_{m_p}$, system response measurements can be measured or computed at predetermined frequency points, $s=\xi_1,\ldots, \xi_{m_s}$.  In the context of our working hypotheses, we will assume that the point sets, $\{\mu_1,\ldots, \mu_{m_p}\}$ and $\{\xi_1,\ldots, \xi_{m_s}\}$, respectively comprise $m_p$ and $m_s$ distinct points in $\mathbb{C}$,  and that the magnitude and phase of the complex frequency response $\cbfH(\xi_i, \mu_j)$ is available for $i=1,\ldots,m_s$ and $j=1,\ldots,m_p$.  These values of $\cbfH(\xi_i, \mu_j)$ will be the only information presumed to be available for the system of interest.  
For projection-based parametric model reduction, which requires access to internal dynamics of the underlying system,  we refer the reader to
\cite{Benner2016,hesthaven2016certified,QuaRM11,benner2017model} and the references therein. 
 Based on the samples $\cbfH(\xi_i, \mu_j)$, we will seek a bivariate function
$\widehat{\cbfH}(s, p)$ such that $\widehat{\cbfH}(\xi_i, \mu_j)\approx\cbfH(\xi_i, \mu_j)$ and such that $\widehat{\cbfH}(s, p)$ itself represents a parameterized transfer function.  We measure closeness in the least-squares sense:
	Given data $\left\{ \xi_i,\mu_j,  \cbfH(\xi_i,\mu_j) \right\}$ with $1\leq i\leq m_s$ and $1\leq j \leq m_p$,
	find a stable bivariate rational function $\widehat\cbfH$  such that 		\begin{equation}
		\sum_{i=1}^{m_s} \sum_{j=1}^{m_p} \left| \widehat\cbfH(\xi_i,\mu_j) - \cbfH(\xi_i,\mu_j) \right|^2 \to \min.
	\end{equation}

 What particular structure should one enforce on $\widehat\cbfH(s, p)$ ? 

\subsection{The solution framework}\label{subsec:Intro_SolnFrmwrk}

We require that $\widehat{\cbfH}(s, p)$ be  a proper rational function with respect to $s$ for each $p$ with poles in the left-half plane, and that for each $s$ at which $\widehat{\cbfH}(s, p)$ is finite,  $\widehat{\cbfH}(s, p)$ have a functional dependence on $p$ that may be polynomial or more general.  

In order to make this more precise, for any integer $n_s\geq 1$, denote by $\calR_{n_s}$ the set of strictly proper, stable, rational functions of order no greater than ${n_s}$  (i.e., for any $g\in\calR_{n_s}$, $g$ is the ratio of polynomials that may be assumed to be relatively prime, with denominator that has polynomial order no greater than ${n_s}$ and roots in the open left half-plane; the numerator has polynomial order strictly less than that of the denominator).  $\calR_{n_s}$ is not itself a vector space, though it is a differentiable manifold containing all subspaces having the form 
$$\mathfrak{R}_{r_s}=\mathsf{span}\{\mathfrak{h}_1(s), \mathfrak{h}_2(s), ..., \mathfrak{h}_{r_s}(s)\}$$
 where $\{\mathfrak{h}_k\}_1^{r_s}$ are proper, stable rational functions with orders summing to a quantity no larger than $n_s$: 
 $\sum_{k=1}^{r_s} \mathsf{order}\left(\mathfrak{h}_k(s)\right)\leq n_s$.
We assume for convenience that $\{\mathfrak{h}_k\}_1^{r_s}$ have mutually distinct poles, so that for some choice, 
$\{\lambda_{\ell}\}_{\ell =1}^{n_s}\subset \mathbb{C}$, distinct points chosen in the open left half-plane, 
$$\mathfrak{R}_{r_s}\subset \mathsf{span}\left\{(s-\lambda_1)^{-1}, (s-\lambda_2)^{-1}, ..., (s-\lambda_{n_s})^{-1}\right\}\subset \calR_{n_s} $$

We allow the parameter $p$ to vary over a convex, compact subset, $\calP\subset \mathbb{C}$, containing our given parameter control values, $\{\mu_1,\ldots,\mu_{m_p}\}\subset \calP$.  Denote by  $\mathcal{F}(\calP)$ the set of continuous functions mapping $\calP$ to $\mathbb{C}$. For ${r_p}\geq 1$ and some choice of functions, 
$\{P_1, P_2, \ldots, P_{r_p} \}\subset \mathcal{F}(\calP)$, let  $\mathfrak{P}_{r_p}=\mathsf{span}\{P_1, P_2, \ldots, P_{r_p}\}\subset \mathcal{F}(\calP)$.  

Consider the tensor product space, $\mathfrak{R}_{r_s} \otimes \mathfrak{P}_{r_p}$, defined formally as the span 
of pairwise elementary products of functions drawn from bases of $\mathfrak{R}_{r_s}$ and $ \mathfrak{P}_{r_p}$:
$\sum_{k\ell}x_{k\ell}\, \mathfrak{h}_k(s)P_{\ell}(p)$, with $1\leq k\leq r_s$, $1\leq \ell \leq r_p$ and $\{x_{k\ell}\}\subset\mathbb{C}$.  
Observe that any function 
\begin{equation} \label{eqn:Hhatform}
\widehat\cbfH(s,p) =\sum_{k=1}^{r_s}\sum_{\ell=1}^{r_p}x_{k\ell}\, \mathfrak{h}_k(s)P_{\ell}(p) \in \mathfrak{R}_{r_s} \otimes \mathfrak{P}_{r_p}
\end{equation}
 will have the properties:
\begin{enumerate}
\item For any fixed $\hat{p}\in \calP$, $\widehat\cbfH(\cdot,\hat{p})$ is a stable transfer function of order $r_s$ (or less) with poles contained in the set $\{\lambda_1, \lambda_2, ..., \lambda_{n_s} \}$.
\item For any fixed $\hat{s}\neq\lambda_{\ell}$, $\ell=1,\ldots, r_s$,  $\widehat\cbfH(\hat{s},\cdot)\in \mathfrak{P}_{r_p}.$ That is, 
$\widehat\cbfH(\hat{s},p)$ has a parametric dependence on $p$ described by $\{P_1, P_2, ..., P_{r_p} \}$.
\end{enumerate}
We begin by considering the problem:
\begin{center}

\begin{minipage}{0.9\textwidth} 
	Given data $\left\{ \xi_i,\mu_j,  \cbfH(\xi_i,\mu_j) \right\}$ with $1\leq i\leq m_s$ and $1\leq j \leq m_p$,
	and subspaces $\mathfrak{R}_{r_s}$ and $\mathfrak{P}_{r_p}$ as described above,  
	find a stable bivariate rational function that fits the data in a least squares sense:\\
	\hspace*{5mm}Find $\widehat\cbfH\in \mathfrak{R}_{r_s} \otimes \mathfrak{P}_{r_p}$ 
	as in \cref{eqn:Hhatform} that solves
	\begin{equation}\label{eq:pvf:bigLSProblem}
		\sum_{i=1}^{m_s} \sum_{j=1}^{m_p} \left| \widehat\cbfH(\xi_i,\mu_j) - \cbfH(\xi_i,\mu_j) \right|^2 \to \min.
	\end{equation}
\end{minipage}
\bigskip
\end{center}
This problem is straightforwardly solved, at least in principle.  Define the matrices 
\begin{equation}  \label{eq:LScoeffMatr}
\begin{array}{c}\mathbb{A}=[\mathfrak{h}_j(\xi_i)]\in \mathbb{C}^{m_s\times r_s},\qquad
\mathbb{B}=[P_j(\mu_i)]\in \mathbb{C}^{m_p\times r_p},\\[2mm]
~~\mbox{and}~~  \mathbb{H}= [\cbfH(\xi_i,\mu_j)]\in\mathbb{C}^{m_s\times m_p}. 
\end{array}
\end{equation}   Then with $X=[x_{ij}]$ and 
$\widehat\cbfH(s,p)=\sum_{k,\ell} x_{k\ell} \,\mathfrak{h}_k(s)P_{\ell}(p)$ as in 
\cref{eqn:Hhatform}
 we observe
$$
\sum_{i=1}^{m_s} \sum_{j=1}^{m_p} \left| \widehat\cbfH(\xi_i,\mu_j) - \cbfH(\xi_i,\mu_j) \right|^2 =
\left\| \mathbb{A}X \mathbb{B}^\top - \mathbb{H}  \right\|_F^2=\|(\mathbb{B}\otimes \mathbb{A})\, \mathsf{vec}(X)-\mathsf{vec}(\mathbb{H})\|_2^2
$$
where we have made use of the Frobenius norm of a matrix: $\|M\|_F^2=\mathsf{trace}(\overline{M}^\top M)$, the stacking operator 
$\mathsf{vec}(\cdot)$, and the matrix Kronecker product, $\otimes$ (see, e.g., Chapter 4 of \cite{HornJohnson1994TopicsMatrAnalysis}). The minimal norm least squares solution to (\ref{eq:pvf:bigLSProblem}) is given by 
$\widehat\cbfH(s,p)=\sum_{k,\ell} \hat{x}_{k\ell} \,\mathfrak{h}_k(s)P_{\ell}(p)$, where $\hat{X}=[\hat{x}_{ij}]$ is obtained from
\begin{equation} \label{eq:FixedBasisLSsoln}
 \mathsf{vec}(\hat{X})=(\mathbb{B}\otimes \mathbb{A})^{\dag}\mathsf{vec}(\mathbb{H})=(\mathbb{B}^{\dag}\otimes \mathbb{A}^{\dag})\mathsf{vec}(\mathbb{H})
\iff \hat{X}= \mathbb{A}^{\dag}\,\cdot\, \mathbb{H} \,\cdot\, (\mathbb{B}^{\dag})^\top.
\end{equation}
where $(\cdot )^\dag$ denotes the matrix pseudo inverse.
This discussion can be easily extended to allow weighting factors $w_{ij}\geq 0$ that assign varying relevance to each data input pair $\xi_i \longleftrightarrow \mu_j$ in the optimization (\ref{eq:pvf:bigLSProblem}). The objective function then becomes 
	\begin{equation}\label{eq:pvf:bigLSProblem-W}
	\sum_{i=1}^{m_s} \sum_{j=1}^{m_p} w_{ij}\left| \widehat\cbfH(\xi_i,\mu_j) - \cbfH(\xi_i,\mu_j) \right|^2 \to \min.
	\end{equation}
Let $W = (\sqrt{w_{ij}})$. Then, using the above notation and the Hadamard (elementwise) product $\circ$, we obtain a weighted least squares problem:
	\begin{eqnarray*}
	\sum_{i=1}^{m_s} \sum_{j=1}^{m_p} w_{ij}\left| \widehat\cbfH(\xi_i,\mu_j) - \cbfH(\xi_i,\mu_j) \right|^2 &=&
	\left\| ( \mathbb{A}X \mathbb{B}^\top - \mathbb{H} ) \circ W  \right\|_F^2 \\ &=& \| \mathrm{diag}(\mathsf{vec}(W)) \left( (\mathbb{B}\otimes \mathbb{A})\, \mathsf{vec}(X)-\mathsf{vec}(\mathbb{H}) \right)\|_2^2 .
	\end{eqnarray*}
Notice in particular that the choice of weight $w_{IJ}=0$ could be appropriate if no observation is available for the case $s=\xi_I$ and $p=\mu_J$ or the data are irrelevant.

One difficulty that has not yet been addressed here lies in representation, that is: (a) what choice to make for frequency response basis functions, $\{\mathfrak{h}_k\}_1^{r_s}$, and more generally, the associated pole locations, $\{\lambda_k\}_{k=1}^{n_s}$, which as yet have only been constrained to lay in the open left half-plane; and (b) what choice to make for parametric bases, $\{P_{\ell}\}_{\ell=1}^{r_p}$, which as yet are only constrained to be continuous functions on $\calP$.   These choices will make up the first of two major concerns in the discussion that follows; our second concern will focus on eliminating (nearly) redundant information that may be implicit in the expression, $\widehat\cbfH(s,p)=\sum_{k,\ell} \hat{x}_{k\ell}\,\mathfrak{h}_k(s)P_{\ell}(p)$, in order to arrive at a more compact representation. 

We propose a two phase approach to address these two difficulties:  the first phase encodes the response data into a high fidelity model of modest order.  To accomplish this, we extend the well-known non-parametric least-squares fitting approaches so as to accommodate parameterized systems. This results in a (discrete) least-squares problem formulated with respect to both frequency and parameter that identifies ``local" system response features.  Our approach produces effective choices for the subspaces $\mathfrak{R}_{r_s}$ and $\mathfrak{P}_{r_p}$ described above.  The second phase that follows constitutes a compression stage that serves to eliminate redundancy in the intermediate model that was obtained in the first phase.   This is accomplished using an 
$\mathbf{\mathcal{H}}_2$-optimal model reduction strategy well-suited for the specialized parametric structure of the intermediate model obtained in the first phase. Our final compressed model inherits the parametric dependence of the intermediate model while maintaining high fidelity, yet our final model will generally have dramatically smaller system order. 

The strategies that we propose have straightforward extensions to the MIMO (multiple-input / multiple-output) setting as well as to multi parameter settings, however for clarity we restrict the discussion here to the SISO - scalar parameter setting, adding a brief discussion of multi parameter case in \cref{sec:pvf:severalParam}.

%

%!TEX root = ms.tex
%  +  +  +  +  +  +  +  +  +  +  +  +  +  +  +  +  +  +  +  +  +  +  +  +  +  +  +  +  +  +  +  +  +  +  +  +  +  +  +  +
\section{Local Models}\label{sec:localmodels}
For any fixed parameter, $\hat{\mu}\in \calP$, a \emph{parametrically localized model} or more briefly, a \emph{local model}, will refer to a transfer function, $\hat{\localH}(s)$, that is intended to model the response of the parameterized system in a neighborhood of $p=\hat{\mu}$.   Within the empirical framework we have adopted here,  a local model, 
$\hat{\localH}(s)$, should match the data well for some choice of $\hat{\mu}\in \{\mu_1,\,\mu_2,\,...,\,\mu_{m_p}\}$: $\hat{\localH}(\xi_i)\approx \cbfH(\xi_i, \hat{\mu})$ for $i=1,...,m_s$.  
One might reasonably expect that a parameterized model $\widehat\cbfH\in \mathfrak{R}_{r_s} \otimes \mathfrak{P}_{r_p}$ that is an effective solution to (\ref{eq:pvf:bigLSProblem}), i.e., a solution that produces a small residual, then  $\widehat\cbfH(s,p)$, will also determine implicitly a family of effective local models as well: $\widehat\cbfH(\cdot,\mu_1),\,\widehat\cbfH(\cdot,\mu_2),...\,\widehat\cbfH(\cdot,\mu_{m_p})$ each should be effective local models.   Since each such local model must also be contained in $\mathfrak{R}_{r_s}$, we may choose to generate an initial subspace $\mathfrak{R}_{r_s}$ by finding an effective local model for each $\{\mu_1,\,\mu_2,\,...,\,\mu_{m_p}\}$.  
Toward that end, for each  $j=1,2,\ldots,m_p$, pick an integer $\nu_j \geq 1$ such that $\sum_j \nu_j \leq n_s$ and define
\begin{equation}\label{eq:pvf:localmodels}
	\localH_j \vcentcolon = \arg\min_{\localH}\sum_{i=1}^{m_s} w_i \left| \localH(\xi_i) - \cbfH(\xi_i,\mu_j) \right|^2,
\end{equation}
where the minimization is taken over strictly proper, stable rational functions having order $\nu_j$ or less.  
As in (\ref{eq:pvf:bigLSProblem-W}), the positive weights $\{w_i\}_1^{m_p}$, provide flexibility in either balancing or focussing the importance of each summand in the objective function, but further also allow one to approximate conditions leading to near-$\mathcal{H}_2$ optimality (see, e.g., \cite{Drmac2015a}).
This process defines a total of $m_p$ local models, so we take $r_s=m_p$ and 
$$\mathfrak{R}_{m_p}=\mathsf{span}\{\mathfrak{h}_1(s), \mathfrak{h}_2(s), ..., \mathfrak{h}_{m_p}(s)\}\subset \calR_{n_s} $$
will be our (initial) frequency-domain factor space that can be used in solving (\ref{eq:pvf:bigLSProblem}).  
We provide an illustration of the full process in Section \ref{subsec:CoupledDataFitting} after discussing how one may obtain effective local models.

\subsection{Local models via Vector Fitting} \label{subsec:LocModviaVF}
The Sanathanan-Koerner iteration \cite{Sanathanan1963b} and its later refinement, Vector Fitting \cite{Approximation1999}, provide effective tools for solving each of the nonlinear least squares problem that are implicit in (\ref{eq:pvf:localmodels}). 
There are other strategies that can be followed as well, e.g., a recent development that is compatible with our solution framework is the \textsf{rkfit} algorithm as described in \cite{berljafa2017rkfit}.  
 We use a Vector~Fitting approach, which we summarize here.   

A strictly proper rational function, $\localH(s)$, of order $\nu\geq 1$ may be represented in barycentric form as
\begin{equation}\label{eq:pvf:ratRepH1}
	\localH(s) = \dfrac{n(s)}{d(s)} = \dfrac{\sum_{k=1}^{\nu} \dfrac{\psi_k}{s-\lambda_k}}{1 + \sum_{k=1}^{\nu} \dfrac{\varphi_k}{s-\lambda_k}}. 
\end{equation}
As (\ref{eq:pvf:ratRepH1}) suggests, $\displaystyle n(s)$ is taken to be the numerator of the rightmost expression and $\displaystyle d(s)$ is the denominator.  The \emph{nodes} $\{\lambda_k\}_1^{\nu}\subset \bbC$ are presumed to be distinct but otherwise are arbitrary, at least in principle.  Evidently, they constitute the poles of both $n(s)$ and $d(s)$ (though \emph{not} of $\localH(s)$). $\psi_k$ and $\varphi_k$ are the associated \emph{residues} of $n(s)$ and $d(s)$, respectively.  The barycentric representation of rational functions has numerical advantages over representations that involve polynomial ratios, \cite{Webb2012,berrut2004barycentric,Higham2004}.  Moreover, the value of $\localH(s)$ at any $\lambda_k$ is directly provided in terms of the residues: $\localH(\lambda_k)=\frac{\psi_k}{\varphi_k}$ for 
$k=1,2,\ldots, \nu$.    For any fixed choice of distinct nodes $\{\lambda_k\}_1^{\nu}$, the minimization in (\ref{eq:pvf:localmodels}) is a \emph{nonlinear least squares problem} with respect to the residues $\{\psi_k\}_1^{\nu}$ and $\{\varphi_k\}_1^{\nu}$.

Observe that the objective function of (\ref{eq:pvf:localmodels}), for a fixed $\mu_j$, can be expressed as
$$
\sum_{i=1}^{m_s} w_i \left| \frac{n(\xi_i)}{d(\xi_i)} - \cbfH(\xi_i,\mu_j) \right|^2=
\sum_{i=1}^{m_s} \frac{w_i}{|d(\xi_i)|^2} \left| n(\xi_i) - \cbfH(\xi_i,\mu_j) d(\xi_i) \right|^2
$$ 
and that the parenthesized quantity is \emph{linear} with respect to the residues $\{\psi_j\}_1^{\nu}$ and $\{\varphi_j\}_1^{\nu}$. 
Following \cite{Sanathanan1963b}, we are led to recast the nonlinear least-squares problem as a sequence of weighted linear least-squares problems:
\begin{equation}\label{eq:pvf:introMinProblem1}
	\sum_{i=1}^{m_s} \dfrac{w_i}{|d^{(k-1)}(\xi_i)|^2} \left|  n^{(k)}(\xi_i) - d^{(k)}(\xi_i)\cbfH(\xi_i,\mu_j) \right|^2\to\min,\qquad k=1,2,\ldots
\end{equation}
which can be formulated in standard form as 
\begin{equation}
	\left\| \bfDelta^{(k-1)} \left( \calA \bfx^{(k)} - \bfb \right) \right\|_2^2 \to \min,\quad k=1,2,\ldots
\end{equation}
\begin{align}
\mbox{where}\quad &
\bfDelta^{(k-1)} = \operatorname{diag}\left\{ \sqrt{w_i}/\left| d^{(k-1)}(\xi_i) \right| \right\}_{i=1}^{m_s} \in\bbC^{m_s\times m_s}, \nonumber\\[2mm]
& \bfx^{(k)} = \begin{bmatrix} \psi_1^{(k)} & \ldots & \psi_{\nu}^{(k)} & \varphi_1^{(k)} & \ldots & \varphi_{\nu}^{(k)} \end{bmatrix}^\top \in \bbC^{2\nu},\quad \mbox{and} \nonumber\\[2mm]
&	\calA \vcentcolon= \begin{bmatrix}
		\dfrac{1}{\xi_1 - \lambda_1} & \cdots & \dfrac{1}{\xi_1 - \lambda_{\nu}} & \dfrac{\cbfH(\xi_1,\mu_j)}{\xi_1 - \lambda_1} & \cdots & \dfrac{\cbfH(\xi_1,\mu_j)}{\xi_1 - \lambda_{\nu}} \\
		\vdots & \ddots & \vdots & \vdots & \ddots & \vdots \\
		\dfrac{1}{\xi_{m_s} - \lambda_1} & \cdots & \dfrac{1}{\xi_{m_s} - \lambda_{\nu}} & \dfrac{\cbfH(\xi_{m_s},\mu_j)}{\xi_{m_s} - \lambda_1} & \cdots & \dfrac{\cbfH(\xi_{m_s},\mu_j)}{\xi_{m_s} - \lambda_{\nu}}
	\end{bmatrix} \in\bbC^{m_s\times 2\nu} \label{eq:SKcalA}
\end{align}\\[2mm]
%%%%%%%%%%%%%%%%%%%%%%%%%%%%%%%%%%%%%%%%%%%%%%%%%%%%%%%%%%%%%%%%%%%%%%%%%%
This describes a step of the \emph{Sanathanan-Koerner iteration}. 

The key innovation of \emph{Vector Fitting} (\textsf{VF}), as developed in \cite{Approximation1999}, that distinguishes it from the Sanathanan-Koerner formulation lies in a clever change of the nodes $\{\lambda_j^{(k)}\}_{j=1}^{\nu}$ at iteration step, $k$, leading to a different representation within each step. Let the matrix $\calA^{(k)}$ be defined analogously to $\calA$  in (\ref{eq:SKcalA}), but with the previously fixed nodes, $\lambda_j$, replaced by $\lambda_j^{(k)}$ that will vary from step to step.  At each step $\localH^{(k)}(s)=n^{(k)}(s)/d^{(k)}(s)$; both $n^{(k)}(s)$ and $d^{(k)}(s)$ are defined using the nodes $\lambda_j^{(k)}$ that are then modified in the following way. 

Consider advancing (\ref{eq:pvf:introMinProblem1}) from step $k$ to step $k+1$.
Since the most recent $d^{(k)}$ is available, compute its zeros, $\{\lambda_j^{(k+1)}\}_{j=1}^{\nu}$, and obtain an equivalent representation:
\begin{equation} \label{PoleZero_d_k}
{d}^{(k)}(s) = 1+  \sum_{j=1}^{\nu} \frac{\varphi_j^{(k)}}{s-\lambda_j^{(k)}}
= \frac{\prod_{j=1}^{\nu} (s-\lambda_j^{(k+1)})}{\prod_{j=1}^{\nu} (s-\lambda_j^{(k)})}.
\end{equation}
This factorization can then be used to rewrite the objective (\ref{eq:pvf:introMinProblem1}) as 
\begin{align}
&\sum_{i=1}^{m_s} \frac{w_i}{|{d}^{(k)}(\xi_i)|^2} \left| \sum_{j=1}^{\nu} \frac{\psi_j^{(k+1)}}{\xi_i-\lambda_j^{(k)}}
- \cbfH(\xi_i) \left(1+\sum_{j=1}^{\nu} \frac{\varphi_j^{(k+1)}}{\xi_i-\lambda_j^{(k)}}\right) \right|^2  \label{eq:2.10}\\
&=\sum_{i=1}^{m_s} w_i \left| \frac{{\prod_{j=1}^r (\xi_i-\lambda_j^{(k)})}}{{\prod_{j=1}^{\nu} (\xi_i-\lambda_j^{(k+1)})}}\right|^2
\left| \frac{\tilde{p}^{(k+1)}(\xi_i)}{{\prod_{j=1}^{\nu} (\xi_i-\lambda_j^{(k)})}} - \cbfH(\xi_i) \frac{\tilde{q}^{(k+1)}(\xi_i)}{{\prod_{j=1}^{\nu} (\xi_i-\lambda_j^{(k)})}}\right|^2  \nonumber
\end{align}
where $\tilde{p}^{(k+1)}$ and $\tilde{q}^{(k+1)}$ are, polynomials of degrees ${\nu}-1$ and ${\nu}$, respectively.
	Continuing with similar algebraic manipulations, one obtains
	\begin{eqnarray}
	\ldots &=&\sum_{i=1}^{m_s} w_i
	\left| \frac{\tilde{p}^{(k+1)}(\xi_i)}{{\prod_{j=1}^{\nu} (\xi_i-\lambda_j^{(k+1)})}} - \cbfH(\xi_i) \frac{\tilde{q}^{(k+1)}(\xi_i)}{{\prod_{j=1}^{\nu} (\xi_i-\lambda_j^{(k+1)})}}\right|^2  \nonumber\\
	&=& \sum_{i=1}^{m_s} w_i \left| \sum_{j=1}^{\nu} \frac{\tilde\psi_j^{(k+1)}}{\xi_i-\lambda_j^{(k+1)}}
	- \cbfH(\xi_i) \left( 1 + \sum_{j=1}^{\nu} \frac{\tilde\varphi_j^{(k+1)}}{\xi_i-\lambda_j^{(k+1)}} \right) \right|^2 \label{tildephi} \\
	&=& \|D_{w}(\calA^{(k+1)} {\tilde{x}^{(k+1)}} - h )\|_2^2\rightarrow \min ,\;\;\;\; D_{w}=\mathrm{diag}(\sqrt{w_i})_{i=1}^{m_s}\nonumber
	\label{eq:VF_unweighted}
	\end{eqnarray}
	where {$\tilde{x}^{(k+1)} = \left( \tilde{\phi}_1^{(k+1)} ~ \tilde{\phi}_2^{(k+1)} ~ \cdots ~ \tilde{\phi}_{\nu}^{(k+1)} ~ \tilde{\varphi}_1^{(k+1)} ~ \tilde{\varphi}_2^{(k+1)} ~ \cdots ~ \tilde{\varphi}_{\nu}^{(k+1)}\right)^\top$}
	with {$\tilde{\phi}_j^{(k+1)}$} and   {$\tilde{\varphi}_j^{(k+1)}$} as defined in \eqref{tildephi} denoting  the 
		coefficients of $\localH^{(k+1)}(s)$ in a barycentric representation having nodes $\lambda_j^{(k+1)}$, $j=1,\ldots, \nu$. 
 \textsf{VF} may be concisely described as a representation of the Sanathanan-Koerner iteration in barycentric form with moving nodes. Upon convergence, the nodes $\lambda_j^{(k)}$ become a fixed point for the iteration, the residues $\varphi_j^{(k)}$ approach zero and at some $k_*$, $\localH^{(k_*)}(s)\approx n^{(k_*)}(s)$ is returned in a convenient  pole-residue form.

%%%%%%%%%%%%%%%%%%%%%%%%%%%%%%%%%%%%%%%%%%%%%%%%%%%%%%%%%%%%%%%%%%%%%%%%%%
More details on implementation and variants of \textsf{VF} can be found in
\cite{Chinea2011,Deschrijver2007a,Deschrijver,Gustavsen2007,Beygi2012,Drmac2015a}.
It is worth noting that convergence of \textsf{VF} in general remains an open problem \cite{Grivet-Talocia2006,Shi2016a,Gustavsen2006,Semlyen2000,Hendrickx2004,Lefteriu2013}.
Nonetheless, in practice we typically observe convergence in relatively few iterations, even with poorly selected initial nodes. 

\subsection{Using Local Models for the Coupled Data Fitting Problem} \label{subsec:CoupledDataFitting}
Once a set of local models, $\{\mathfrak{h}_1(s), \mathfrak{h}_2(s), ..., \mathfrak{h}_{m_p}(s)\}$, has been obtained 
associated with the given parameter sampling, $\{\mu_1,\,\mu_2,\,...,\,\mu_{m_p}\}$,
these local models can be used to define a factor space, 
$\mathfrak{R}_{m_p}=\mathsf{span}\{\mathfrak{h}_1(s), \mathfrak{h}_2(s), ..., \mathfrak{h}_{m_p}(s)\}$, that may be 
expected to describe the range of variability of the data with respect to frequency.   In order to solve (\ref{eq:pvf:bigLSProblem}), 
we need also to posit a parameter dependence through the specification of a parameter factor space, 
$\mathfrak{P}_{r_p}$, that serves to knit together the local models that we have developed in \S\ref{sec:localmodels}.
The simplest model of parameter dependence involves choosing an \emph{a priori} fixed subspace spanned 
by polynomials or another fixed basis: $\mathfrak{P}_{r_p}=\mathsf{span}\{P_1, P_2, \ldots, P_{r_p}\}$.  
This is discussed in \S\ref{subsec:Intro_SolnFrmwrk}; we summarize the aggregated strategy as 
\cref{alg:pvf:paramVF1} for the unweighted case ($w_{ij}=1$).  
\begin{algorithm}[H]
	\caption{Parametric Fitting - Fixed Parametric Basis}\label{alg:pvf:paramVF1}
		\begin{tabularx}{\textwidth}{lX}
		{\bf INPUT:}&  Observed response data: $\left\{ \xi_i,\mu_j,\cbfH(\xi_i,\mu_j)\right\}_{i=1,j=1}^{i=m_s,j=m_p}$, \\
		& Target model order: $n_s$\\  & Parametric basis functions $\left\{ P_\ell \right\}_{\ell=1}^{r_p}$ \\[2mm]
		{\bf OUTPUT: } & Bivariate function $\widehat\cbfH(s,p)$ that fits the data and solves (\ref{eq:pvf:bigLSProblem}).
	\end{tabularx}
	\begin{enumerate}
		\item Choose \emph{local model orders}, $\nu_j\geq 1$ such that $\sum_{j=1}^{m_p} \nu_j =n_s$. 
		\item Construct \emph{local models}, $\localH_j(s)$, of order $\nu_j$ for each parameter $\mu_j$ following (\ref{eq:pvf:localmodels}).
		\item Solve the least squares problem (\ref{eq:pvf:bigLSProblem}) using (\ref{eq:FixedBasisLSsoln}) to obtain the coefficients, $\hat{X}=[\hat{x}_{ij}]$.
	        \item Return the parametrized intermediate model:  given by 
$$	\widehat\cbfH(s,p) = \sum_{k=1}^{m_p} \sum_{\ell=1}^{r_p} \hat{x}_{k\ell} P_\ell(p) \localH_k(s). $$
	\end{enumerate}
\end{algorithm} 
 Combining local reduced models with various parametric bases $\mathfrak{P}_{r_p}$ has been considered in \cite{morBauB09,baur2011parameter}. However, these works focused on interpolatory bases, i.e., $\widehat{\cbfH}(s,p)$ interpolated the local model  $\mathfrak{h}_k(s)$  at the parameter sample $p = \mu_k$, for $k=1,2,\ldots,m_p$. Here we focus on parametric bases to minimize a discrete joint least-squares measure. And more importantly, in \Cref{sec:pvf:paramVF}, we will  allow these bases functions in $p$ to adaptively vary. 
For the concept of  local reduced models in a projection-based setting, we refer the reader to, e.g., \cite{Benner2016,Amsallem2008,Degroote2010,Lohmann2009,Panzer_etal2010} and the references therein.
 \begin{example}   \label{ex:pvf:polBeamExample}
 We illustrate \cref{alg:pvf:paramVF1} on data generated from an idealized vibration model of a cantilevered Timoshenko beam \cite{panzer2009generating} with proportional damping: \begin{equation}\label{eq:pvf:2ndOrdEx}
	\cbfH(s,p)= \bfc^\top\left(\bfM s^2 + \left({\textstyle\frac12}\bfM + p \bfK \right) s + \bfK\right)^{-1} \bfb.
\end{equation}
The matrices $\bfM$ and $\bfK$ represent the distributed mass and stiffness of a Timoshenko beam approximated with finite elements; 
for the model considered here, we have 2300 degrees of freedom. 
The parameter $p$ represents a  damping parameter and we consider $p\in \calP=[0,0.8]$.
 We sample at $m_s = 80$ logarithmically spaced frequency points between $10^{-3}$ and $10^3$ on the imaginary axis, capturing the main dynamic range of the model and consider $m_p=10$ parameter sampling points, equally spaced across the interval $[0.01,0.8]$.

 Local models are constructed at each of the parameter sampling points using \textsf{VF} as described above. Each local model has order $\nu_j=10$  
 producing an aggregate rational model, $\widehat\cbfH(s,p)$, having $s$-order  $n_s=100$. 

We take for $\mathfrak{P}_{r_p}$ the subspace of polynomials having order up to 5  (so that $r_p=6$); 
we use Bernstein polynomials as basis elements, $P_\ell(p)$; see \cite{Balazs1982,Farouki2003}.

In \cref{eq:pvf:p1DBeamErrorPoly}, we show the Bode amplitude and error plots at two representative parameter values, $p=0.01$ and $p=0.22$.
The original sampled function is shown in blue, our approximation is shown with dashed red lines, and (absolute) pointwise error appears in green.
To illustrate the quality of the parametric reduced model over the whole parameter domain $p\in \calP=[0,1]$,  in \cref{fig:pvf:p1DBeamErrorPolyErrCont}, we show the relative error between $\widehat\cbfH(s,p)$ and the original sampled model $\cbfH(s,p)$ with respect to the $\cbfH_2$ and $\cbfH_\infty$ norm for every parameter value, i.e., for  $\hat{p}\in \calP=[0,1]$, we plot

{\small
\begin{equation*} 
\begin{aligned}
\left\| \cbfH(\cdot,\hat{p}) - \widehat\cbfH(\cdot,\hat{p}) \right\|_{\cbfH_2} 
& = \left(\frac{1}{2\pi} \int_{-\infty}^\infty \left| \cbfH(\rmi\omega,\hat{p}) - \widehat\cbfH(\rmi\omega,\hat{p}) \right|^2 \rmd \omega\right)^{1/2}, \quad \mbox{and}\\[2mm]
\left\| \cbfH(\cdot,\hat{p}) - \widehat\cbfH(\cdot,\hat{p}) \right\|_{\cbfH_\infty} 
& = \sup_{\omega \in \mathbb{R}} \left| \cbfH(\rmi\omega,\hat{p}) - \widehat\cbfH(\rmi\omega,\hat{p}) \right|.
\end{aligned}
\end{equation*} }
We observe that even away from the neighborhood of sampling points, 
the approximation performs quite well. Indeed, over almost the entire parameter domain, we obtain relative accuracy in both measures below $10^{-3}$, with the largest relative error being approximately 
$5\times 10^{-3}$.
	\begin{figure}[H]
		\centering
		\begin{minipage}{0.49\textwidth}
			\ifuseTikzGraphs
			\begin{adjustbox}{max width=0.9\textwidth,center}
				\pgfplotsset{ymin=5e-7, ymax=1e+2, xlabel={Frequency $\omega$},ylabel={Magnitude}}
				\input{"Graphics/PVF_1P/BeamRev_rs10_Poly_NS1.tex"} 
			\end{adjustbox}
			\else\includegraphics[width=0.9\textwidth, trim= 0 0 0 22, clip]{PVF_1P/BeamRev_rs10_Poly_NS1}\fi
			\subcaption{$p=0.01$}
		\end{minipage}\hfill
		\begin{minipage}{0.49\textwidth}
			\ifuseTikzGraphs
			\begin{adjustbox}{max width=0.9\textwidth,center}
				\pgfplotsset{ymin=5e-7, ymax=1e+2, xlabel={Frequency $\omega$},ylabel={Magnitude}}
				\input{"Graphics/PVF_1P/BeamRev_rs10_Poly_NS2.tex"}
			\end{adjustbox}
			\else\includegraphics[width=0.9\textwidth, trim= 0 0 0 22, clip]{PVF_1P/BeamRev_rs10_Poly_NS2}\fi
			\subcaption{$p=0.22$}
		\end{minipage}
		\caption{Frequency responses of $\cbfH(s,p)$ (\textcolor{blue}{\rule[0.025in]{0.12in}{1.5pt}}), $\widehat{\cbfH}(s,p)$ (\textcolor{red}{\rule[0.025in]{0.05in}{1.5pt}\,\rule[0.025in]{0.05in}{1.5pt}}), and the error function
		$\cbfH(s,p)- \widehat{\cbfH}(s,p)$ (\textcolor{green}{\rule[0.025in]{0.12in}{1.5pt}}) at $p=0.01$ and $p=0.22$}
		\label{eq:pvf:p1DBeamErrorPoly}
	\end{figure}
	\begin{figure}[H]
		\centering
		\ifuseTikzGraphs
			\begin{adjustbox}{max width=0.7\textwidth}
\definecolor{mycolor1}{rgb}{0.00000,0.44700,0.74100}%
\definecolor{mycolor2}{rgb}{0.85000,0.32500,0.09800}%
\begin{tikzpicture}

\begin{axis}[%
scale only axis,
xmin=0,
xmax=0.8,
xlabel style={font=\color{white!15!black}},
xlabel={Parameter},
ymode=log,
ymin=0.0001,
ymax=0.00243725675312876,
yminorticks=true,
ylabel style={font=\color{white!15!black}},
ylabel={(local) relative Error},
axis background/.style={fill=none},
title style={font=\bfseries},
title={Approximation Quality over Parameter Range for Polynomial $P_\ell$},
xmajorgrids,
ymajorgrids,
yminorgrids,
legend style={legend cell align=left, align=left, draw=white!15!black}
]
\addplot [color=mycolor1, line width=2.0pt]
  table[row sep=crcr]{%
0	0.000618890508268892\\
0.00808080808080808	0.000324555571321873\\
0.0161616161616162	0.0002106627411995\\
0.0242424242424242	0.000181466000581618\\
0.0323232323232323	0.000199416483443648\\
0.0404040404040404	0.000230715337355694\\
0.0484848484848485	0.000260722649326538\\
0.0565656565656566	0.000285203517491966\\
0.0646464646464646	0.000303339342169165\\
0.0727272727272727	0.000315358965218968\\
0.0808080808080808	0.00032180785288987\\
0.0888888888888889	0.000323335722435294\\
0.096969696969697	0.000320577841892704\\
0.105050505050505	0.000314157864558205\\
0.113131313131313	0.000304657858300256\\
0.121212121212121	0.000292651880159156\\
0.129292929292929	0.000278651383376691\\
0.137373737373737	0.000263162364599327\\
0.145454545454545	0.000246656357658377\\
0.153535353535354	0.000229582967921488\\
0.161616161616162	0.000212380168805975\\
0.16969696969697	0.000195475583423278\\
0.177777777777778	0.000179276531529913\\
0.185858585858586	0.000164190456144493\\
0.193939393939394	0.000150601306972962\\
0.202020202020202	0.000138863348837929\\
0.21010101010101	0.000129264685932491\\
0.218181818181818	0.000121983073317294\\
0.226262626262626	0.00011704542963613\\
0.234343434343434	0.00011429382747041\\
0.242424242424242	0.000113409606124123\\
0.250505050505051	0.000113973992135356\\
0.258585858585859	0.000115542632464562\\
0.266666666666667	0.000117713180616708\\
0.274747474747475	0.000120140374488343\\
0.282828282828283	0.000122570783621407\\
0.290909090909091	0.000124807621265228\\
0.298989898989899	0.000126780666635681\\
0.307070707070707	0.000128406668299375\\
0.315151515151515	0.000129730673797982\\
0.323232323232323	0.000130784358750818\\
0.331313131313131	0.000131640081342783\\
0.339393939393939	0.000132403950296822\\
0.347474747474747	0.000133189804726858\\
0.355555555555556	0.000134097195909105\\
0.363636363636364	0.000135246321441794\\
0.371717171717172	0.000136723097475017\\
0.37979797979798	0.000138587060612422\\
0.387878787878788	0.000140866004144226\\
0.395959595959596	0.000143559960157812\\
0.404040404040404	0.000146630862813038\\
0.412121212121212	0.000150016085499815\\
0.42020202020202	0.000153618498208159\\
0.428282828282828	0.000157328820670239\\
0.436363636363636	0.000161023795122155\\
0.444444444444444	0.000164578748800994\\
0.452525252525253	0.000167864975897671\\
0.460606060606061	0.000170762411943093\\
0.468686868686869	0.000173165589581429\\
0.476767676767677	0.000174973429299903\\
0.484848484848485	0.000176108892735148\\
0.492929292929293	0.000176521141792149\\
0.501010101010101	0.000176172619131526\\
0.509090909090909	0.000175052988720194\\
0.517171717171717	0.00017319932433456\\
0.525252525252525	0.000170662532167938\\
0.533333333333333	0.000167546820610108\\
0.541414141414141	0.000163999330789591\\
0.54949494949495	0.000160210659046465\\
0.557575757575758	0.000156425567378946\\
0.565656565656566	0.000152933223193115\\
0.573737373737374	0.000150056055109335\\
0.581818181818182	0.000148132792008075\\
0.58989898989899	0.000147473231197148\\
0.597979797979798	0.000148319346355869\\
0.606060606060606	0.000150795691471602\\
0.614141414141414	0.00015490879117057\\
0.622222222222222	0.000160482291795383\\
0.63030303030303	0.000167227781516254\\
0.638383838383838	0.000174766231757625\\
0.646464646464647	0.000182641022213237\\
0.654545454545455	0.000190396893014577\\
0.662626262626263	0.0001975402271406\\
0.670707070707071	0.000203618095861675\\
0.678787878787879	0.000208182667456353\\
0.686868686868687	0.000210890477906438\\
0.694949494949495	0.000211496702919545\\
0.703030303030303	0.000209956505902251\\
0.711111111111111	0.000206527525137194\\
0.719191919191919	0.000202062512852626\\
0.727272727272727	0.000198244595678643\\
0.735353535353535	0.000198024461079648\\
0.743434343434344	0.000205613182605787\\
0.751515151515152	0.000225683674261972\\
0.75959595959596	0.0002618237542337\\
0.767676767676768	0.000315652227367075\\
0.775757575757576	0.000387413110924169\\
0.783838383838384	0.000476979041951128\\
0.791919191919192	0.000584409245374149\\
0.8	0.000710107398589854\\
};
\addlegendentry{$\cbfH_2$ - Norm}

\addplot [color=mycolor2, line width=2.0pt]
  table[row sep=crcr]{%
0	0.00243725675312876\\
0.00808080808080808	0.000720053896931061\\
0.0161616161616162	0.000391620106137788\\
0.0242424242424242	0.00079865865274357\\
0.0323232323232323	0.00101990820896652\\
0.0404040404040404	0.00112782598795641\\
0.0484848484848485	0.00116446672658287\\
0.0565656565656566	0.00115524132465305\\
0.0646464646464646	0.00111623962246731\\
0.0727272727272727	0.00105803717858641\\
0.0808080808080808	0.000987778428238408\\
0.0888888888888889	0.000910395115522357\\
0.096969696969697	0.000829341938702432\\
0.105050505050505	0.000748167077454054\\
0.113131313131313	0.000684920478955284\\
0.121212121212121	0.000622401444976451\\
0.129292929292929	0.00055810318404294\\
0.137373737373737	0.000493338119368306\\
0.145454545454545	0.000429203959769293\\
0.153535353535354	0.000367634776085946\\
0.161616161616162	0.000314285681134164\\
0.16969696969697	0.000262990032555574\\
0.177777777777778	0.000218131264813641\\
0.185858585858586	0.000187620935059653\\
0.193939393939394	0.000180742920513518\\
0.202020202020202	0.000175399044405219\\
0.21010101010101	0.000173707935413348\\
0.218181818181818	0.000172771812472927\\
0.226262626262626	0.000170692246596335\\
0.234343434343434	0.000187766367632584\\
0.242424242424242	0.000210859607995084\\
0.250505050505051	0.000230512029117062\\
0.258585858585859	0.000246756077981709\\
0.266666666666667	0.0002596895587301\\
0.274747474747475	0.000271371821523076\\
0.282828282828283	0.000280905897201809\\
0.290909090909091	0.000287660745746576\\
0.298989898989899	0.000291815334329171\\
0.307070707070707	0.000293552723633533\\
0.315151515151515	0.00029307380717724\\
0.323232323232323	0.000290579189172438\\
0.331313131313131	0.000286276515689877\\
0.339393939393939	0.000283113998371609\\
0.347474747474747	0.000279868722979256\\
0.355555555555556	0.000275490692416309\\
0.363636363636364	0.000270153940234051\\
0.371717171717172	0.000264028748749272\\
0.37979797979798	0.000258139622413392\\
0.387878787878788	0.0002545625351593\\
0.395959595959596	0.000250491483536865\\
0.404040404040404	0.000246009086826183\\
0.412121212121212	0.000241193488212123\\
0.42020202020202	0.0002361012732478\\
0.428282828282828	0.00023078184376655\\
0.436363636363636	0.000227892062112095\\
0.444444444444444	0.000236274039599845\\
0.452525252525253	0.000244116959705134\\
0.460606060606061	0.000251049282580028\\
0.468686868686869	0.000256639683975752\\
0.476767676767677	0.000260598179951687\\
0.484848484848485	0.00026263565601489\\
0.492929292929293	0.000262867227078926\\
0.501010101010101	0.000261032883321856\\
0.509090909090909	0.000257087632351848\\
0.517171717171717	0.000251002022594072\\
0.525252525252525	0.000242772674737682\\
0.533333333333333	0.000232401837536337\\
0.541414141414141	0.000219912727662546\\
0.54949494949495	0.000205349171829001\\
0.557575757575758	0.000188901075645612\\
0.565656565656566	0.000172665124356867\\
0.573737373737374	0.000158130567389227\\
0.581818181818182	0.000154482682420032\\
0.58989898989899	0.000156891139315127\\
0.597979797979798	0.000160432060884731\\
0.606060606060606	0.000165111605568808\\
0.614141414141414	0.000170808866325554\\
0.622222222222222	0.000177174838227768\\
0.63030303030303	0.000184075833166678\\
0.638383838383838	0.00019100082848325\\
0.646464646464647	0.000197569147247164\\
0.654545454545455	0.000203374003537201\\
0.662626262626263	0.000207974043119444\\
0.670707070707071	0.000210815126735273\\
0.678787878787879	0.000216379210306507\\
0.686868686868687	0.000241734769271575\\
0.694949494949495	0.000264948029073377\\
0.703030303030303	0.000285528002793495\\
0.711111111111111	0.0003029614500563\\
0.719191919191919	0.000316697029255113\\
0.727272727272727	0.000327298811290635\\
0.735353535353535	0.000333700786849677\\
0.743434343434344	0.000338044322989702\\
0.751515151515152	0.000343125571780127\\
0.75959595959596	0.000353238216491484\\
0.767676767676768	0.000369149984689289\\
0.775757575757576	0.000407009518630983\\
0.783838383838384	0.000511544959146191\\
0.791919191919192	0.000633943863635382\\
0.8	0.000774692819107599\\
};
\addlegendentry{$\cbfH_\infty$ - Norm}

\end{axis}
\end{tikzpicture}%
			\end{adjustbox}
		\else\includegraphics[width=0.9\textwidth]{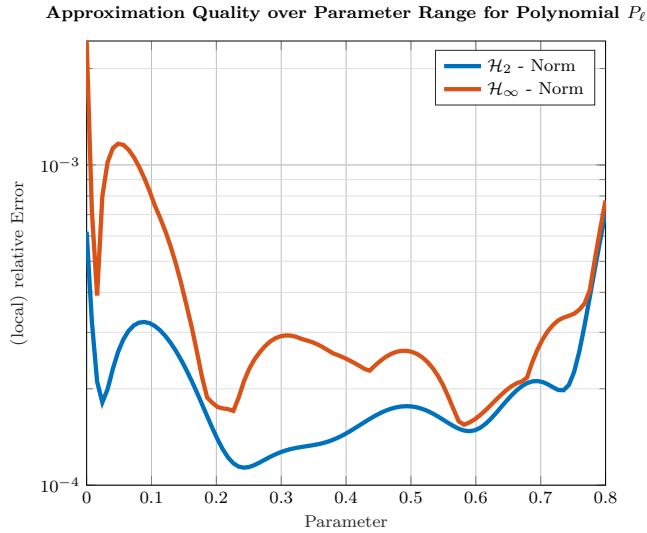}\fi
		\caption{Relative $\cbfH_2$ / $\cbfH_\infty$ errors over the parameter range $[0.01,0.8]$ for parametric fitting using polynomial functions $P_\ell(p)$. Here $r_p = 10$.  }
		\label{fig:pvf:p1DBeamErrorPolyErrCont}
	\end{figure}
	\end{example}

%
%

%!TEX root = ms.tex
%
%
\section{Coupled Parametric Fitting}\label{sec:pvf:paramVF}

We have formulated a solution strategy to the underlying least squares data fitting problem (\ref{eq:pvf:bigLSProblem}) that assumes the form of a solution to be given as a linear combination of products of optimized local models, as defined in (\ref{eq:pvf:localmodels}), with a fixed set of parametric basis functions, $\{P_1, P_2, \ldots, P_{r_p} \}$, that nominally have been defined \emph{independently} of the data. This solution strategy has been formalized as \cref{alg:pvf:paramVF1}.  In this section, we consider an important further refinement that allows the parametric basis functions also to be \emph{adapted} (iteratively) to the data. 

\subsection{Parametric Bases via Variable Projection}
Assume that the family of parametric functions, $\{P_1(p), P_2(p), \ldots, P_{r_p}(p) \}$ depend smoothly on a vector of ancillary (meta-)parameters, $\bfpi\in \mathcal{D}$ for some bounded domain $\mathcal{D}\subset\mathbb{C}^{r_p}$; 
we denote this dependence as $\{P_1^\bfpi(p), P_2^\bfpi(p), \ldots, P_{r_p}^\bfpi(p) \}$. Our focus will be on the case that $\{P_j^\bfpi\}$ are rational functions with pole locations given by the components of $\bfpi\in\mathcal{D}$, but more general settings can be considered, so we keep our discussion general for the time being.  Note that the factor matrix, $\mathbb{B}$,  defined in (\ref{eq:LScoeffMatr}) will now depend smoothly on $\bfpi$ as well: $\mathbb{B}(\bfpi)=[P_j^\bfpi(\mu_i)]\in \mathbb{C}^{m_p\times r_p}$.   We suppose that 
$\mathsf{rank}[\mathbb{B}(\bfpi)]=\dim\mathfrak{P}_{r_p}^\bfpi$ is constant for $\bfpi\in\mathcal{D}$. The pseudoinverse, $\mathbb{B}(\bfpi)^\dag$ and the orthogonal projection onto $\mathsf{Ran}[\mathbb{B}(\bfpi)]$, $\mathbb{Q}(\bfpi)=\mathbb{B}(\bfpi)\mathbb{B}(\bfpi)^{\dag}$,  also vary smoothly with respect to $\bfpi\in\mathcal{D} $, as a consequence. 

The data fitting problem (\ref{eq:pvf:bigLSProblem}) can now be reformulated as:
{\small  \begin{center}
\begin{minipage}{0.9\textwidth} 
Find  $\widehat\cbfH(s,p)=\sum_{k,\ell} \hat{x}_{k\ell} \,\mathfrak{h}_k(s)P_{\ell}^{\hat{\bfpi}}(p)$, \\[2mm]
\hspace*{10mm}where $\hat{X}=[\hat{x}_{ij}]$ and 
 $\hat{\bfpi}\in \mathcal{D}$ solves
\begin{equation} \label{eq:AdaptiveBasisLSsoln1}
 (\hat{X},\, \hat{\bfpi}) =\argmin_{X,\bfpi} \left\| \mathbb{A}\,X\, \mathbb{B}(\bfpi)^\top - \mathbb{H}\right\|_F^2.
\end{equation}
\end{minipage}
\bigskip
\end{center} }

This is a \emph{separable least squares problem} that is linear with respect to $X$ and nonlinear with respect to $\bfpi$.   The linear variables, $X$, can be eliminated, leaving an equivalent nonlinear optimization problem of greatly reduced dimension expressed solely with respect to $\bfpi$.  
Writing $\mathbb{P}=\mathbb{A}\mathbb{A}^{\dag}$ for the orthogonal projection onto $\mathsf{Ran}[\mathbb{A}]$,  we have the equivalent problem:
\begin{equation} \label{eq:AdaptiveBasisLSsoln2}
\begin{array}{c}
\hat{X}= \mathbb{A}^{\dag}\,\cdot\, \mathbb{H} \,\cdot\, (\mathbb{B}(\hat{\bfpi})^{\dag})^\top, \quad \mbox{where}\\[2mm]
\begin{array}{l}
\hat{\bfpi} =\argmin_{\bfpi\in\mathcal{D}}\left\|\mathbb{P}\,\cdot\, \mathbb{H} \,\cdot\, \mathbb{Q}(\bfpi) - \mathbb{H}\right\|_F^2\\[2mm]
     \hspace*{5mm} = \argmin_{\bfpi\in\mathcal{D}}\left\|\left(\mathbf{I}-(\mathbb{Q}(\bfpi)\otimes \mathbb{P})\right)\, \mathsf{vec}(H)\right\|_2^2 .
\end{array}
\end{array}
\end{equation}
Note that $\mathbb{Q}(\bfpi)\otimes \mathbb{P}$ is an orthogonal projector that is dependent on $\bfpi$, and so,  the residual vector, $\mathbf{r}(\bfpi)=\left(\mathbf{I}-(\mathbb{Q}(\bfpi)\otimes \mathbb{P})\right)\, \mathsf{vec}(H)$, is the orthogonal projection of the data, $\mathsf{vec}(H)$, onto the orthogonal complement of the tensor product space $\mathsf{Ran}[\mathbb{B}(\bfpi)] \otimes \mathsf{Ran}[\mathbb{A}]$. 
The map $\bfpi\mapsto \mathbf{r}(\bfpi)$ is a \emph{nonlinear}, albeit smooth, map and  
the minimization appearing in (\ref{eq:AdaptiveBasisLSsoln2}) is thus a nonlinear least squares problem typical of the method of variable projection developed by Golub and Pereyra \cite{golubpereyra1973diffPseudoinv}.   The key observation of Golub and Pereyra was that many nonlinear least squares problems were in fact linear with respect to a substantial number of variables and nonlinear only with respect to a comparative few.  The resulting method of 
 variable projection stemming from this has proven to be quite useful in a wide variety of contexts and there have been substantial subsequent refinements; see e.g., 
\cite{Pereyra2002,Chung2010,olearyrust2013varpro,hokanson2018data}.   
 
 This leads us to our second algorithm which can be viewed simply as a refinement of \cref{alg:pvf:paramVF1}: 
 \begin{algorithm}[H]
 	\caption{Parametric Fitting - Adaptive Parametric Basis}\label{alg:pvf:paramVFVarProj}
		\begin{tabularx}{\textwidth}{lX}
		{\bf INPUT:}&  Observed response data: $\left\{ \xi_i,\mu_j,\cbfH(\xi_i,\mu_j)\right\}_{i=1,j=1}^{i=m_s,j=m_p}$, \\
		& Target model order: $n_s$\\  & Parametric family of $\bfpi$-dependent basis functions $\left\{ P_\ell^{\bfpi} \right\}_{\ell=1}^{r_p}$ \\[2mm]
		{\bf OUTPUT: } & Bivariate function $\widehat\cbfH(s,p)$ that fits the data and solves (\ref{eq:pvf:bigLSProblem}).
	\end{tabularx}
	\begin{enumerate}
		\item Choose \emph{local model orders}, $\nu_j\geq 1$ such that $\sum_{j=1}^{m_p} \nu_j =n_s$. 
		\item Construct \emph{local models}, $\localH_j(s)$, of order $\nu_j$ for each parameter $\mu_j$ following (\ref{eq:pvf:localmodels}).
	         \item Using the \emph{method of variable projection} (see \cite{Chung2010,olearyrust2013varpro}) find 
	            $\hat{X}=[\hat{x}_{ij}]$ and $\hat{\bfpi}$ that solve (\ref{eq:AdaptiveBasisLSsoln1}). 
	        \item Return the parametrized intermediate model given by
	  $\widehat\cbfH(s,p)=\sum_{k,\ell} \hat{x}_{k\ell} \,\mathfrak{h}_k(s)P_{\ell}^{\hat{\bfpi}}(p)$
    \end{enumerate}
\end{algorithm} 

\subsection{Adaptive Rational Parametric Bases}
One simple and often very effective choice for a parametric family of $\bfpi$-dependent basis functions is to choose  $P_\ell^{\bfpi}(p)$ to be 
rational functions with respect $p$, with pole locations given by the components of $\bfpi\in\mathcal{D}$.  We consider this in more detail here and assume for simplicity that 
$$
\bfpi= \begin{bmatrix} \pi_1 & \cdots & \pi_{r_p} \end{bmatrix}^\top \in\mathcal{D}\subset  \bbC^{r_p}\quad\mbox{and}\quad P_\ell^{\bfpi}(p)=\frac{1}{p-\pi_\ell}\quad\mbox{for}\quad \ell=1,\,...,\,r_p.
$$   Assume further that $\mathcal{D}\cap\mathcal{P}=\emptyset$, so that all functions $P_\ell^{\bfpi}(p)$ are smoothly varying with respect to $p$ throughout the parameter range $\calP$.  For any $\bfpi= \begin{bmatrix} \pi_1 & \cdots & \pi_{r_p} \end{bmatrix}^\top \in\mathcal{D}$, using the definition 
$\mathbb{B}$ in \cref{eq:LScoeffMatr}, we can define the factor matrix, 
$$
\mathbb{B}(\bfpi)\vcentcolon =\left[\frac{1}{\mu_i-\pi_j}\right]\in \mathbb{C}^{m_p\times r_p}.
$$

In the following three examples, we illustrate the performance of \cref{alg:pvf:paramVFVarProj} using rational parametric bases as described above. 
In these and subsequent examples, we use an implementation of Variable Projection for Step 3 of \cref{alg:pvf:paramVFVarProj} following \cite{Chung2010}, incorporating a Gauss-Newton strategy to solve the (reduced) nonlinear least squares problem.

 %%%==========================================================================================
\begin{example}\label{ex:pvf:synthRational}
We construct a parametric  transfer function with rational-dependence in $p$:
	\begin{equation} \label{eq:Hratinp}
		\cbfH(s,p) = \sum_{k=1}^{6} \dfrac{\phi_k}{p-\pi_k} G_k(s),\qquad \pi_k\in\bbC,
	\end{equation}
	where the parameter poles are chosen as  $\pi_1 = 0.4$, $\pi_{2,3} =  2\pm 1.5\rmi $, $\pi_{4,5}  = 4\pm 0.8\rmi $, and 
	$\pi_6 = 5.1$.
The  non-parametric rational functions $G_k(s)$ in \cref{eq:Hratinp}  are based on a variation of the Penzl's example \cite[Ex. 3]{morPen99}. Let
\begin{equation}
\begin{aligned}
	\bfA_1 \vcentcolon &=  \begin{bmatrix} -1 & 100 \\ -100 & -1  \end{bmatrix},\quad
	\bfA_2 \vcentcolon =  \begin{bmatrix} -1 & 200 \\ -200 & -1  \end{bmatrix}, \quad
	\bfA_3 \vcentcolon =  \begin{bmatrix} -1 & 400 \\ -400 & -1  \end{bmatrix}, \\
	\bfA_4 \vcentcolon &=  \operatorname{diag}\begin{bmatrix} -1 & -2 & \cdots & -20 \end{bmatrix}, \quad
	\bfA (\zeta) \vcentcolon=\operatorname{diag}\begin{bmatrix}  (\zeta+1)^2\bfA_1 &  \bfA_2 & \bfA_3 & \zeta \bfA_4 \end{bmatrix}, \\
	\bfb(\zeta) \vcentcolon &= \begin{bmatrix} 10 & 10 & \cdots & 10 & 0 & \cdots & \zeta+1 \end{bmatrix},\qquad \mbox{and}\qquad \bfc(\zeta) \vcentcolon =\bfb(\zeta)^\top,
\end{aligned}
\end{equation}
and define the transfer function
\begin{equation}
	 \cbfG(s,\zeta) \vcentcolon = \bfc(\zeta) ^\top \left(s \bfE - \bfA(\zeta) \right)^{-1} \bfb(\zeta). 
\end{equation}
We evaluate  $\cbfG(s,\zeta) $ at $6$ fixed $\zeta_k$ values  
to construct $G_k(s)$ in \cref{eq:Hratinp}: 
\[
G_k(s)  \vcentcolon =  \cbfG(s,\zeta_k), \quad \mbox{where} \quad 
	\begin{bmatrix}
 \zeta_1 & \zeta_2 & \cdots & \zeta_6	\end{bmatrix} =
	\begin{bmatrix}
		0 &  0.29  &  0.57 &   0.86 &   1.14  &  1.43   & 1.71   & 2
	\end{bmatrix}.
\]
We  sample $\cbfH(s,p)$ in \cref{eq:Hratinp} at $100$ frequency samples, logarithmically spaced between $10^{-1}$ and $10^{5}$ on the imaginary axis and $8$ parameter samples linearly spaced on $\calP=[1,5]$. We apply both the fixed polynomial basis (\Cref{alg:pvf:paramVF1}) and adaptive rational basis (\Cref{alg:pvf:paramVFVarProj}) approaches  to $\cbfH(s,p)$. \Cref{fig:pvf:synthPolApprox} and \cref{fig:pvf:synthRatApprox}
show  the amplitude Bode plot comparisons at two of the sampled  points of $p=1.57$ and $p=3.86$, respectively, where the green line represents the error function. 

		\begin{figure}[H]
		\begin{minipage}{0.49\textwidth}
			\ifuseTikzGraphs
				\begin{adjustbox}{max width=0.9\textwidth}
				\pgfplotsset{ymin=0.001, ymax=10, xlabel={Frequency $\omega$},ylabel={Magnitude},title={}}
				\input{"Graphics/PVF_SynthP/PolyApprox_r4_sPt2.tex"}
			\end{adjustbox}
			\else\includegraphics[width=0.9\textwidth, trim= 0 0 0 0, clip]{PVF_SynthP/PolyApprox_r4_sPt2}\fi
			\subcaption{$p=1.57$ -- \Cref{alg:pvf:paramVF1}}
		\end{minipage}\hfill
				\begin{minipage}{0.49\textwidth}
			\ifuseTikzGraphs
				\begin{adjustbox}{max width=0.9\textwidth}
				\pgfplotsset{ymin=0.00031, ymax=10, xlabel={Frequency $\omega$},ylabel={Magnitude}} % 
				\input{"Graphics/PVF_SynthP/RatApprox_r4_sPt2.tex"}
				\end{adjustbox}
			\else\includegraphics[width=0.9\textwidth, trim= 0 0 0 0, clip]{PVF_SynthP/RatApprox_r4_sPt2}\fi
			\subcaption{$p=1.57$ -- \Cref{alg:pvf:paramVFVarProj}}
		\end{minipage}
				\caption{Frequency comparisons of $\cbfH(s,p)$ (\textcolor{blue}{\rule[0.025in]{0.12in}{1.5pt}}), $\widehat{\cbfH}(s,p)$ (\textcolor{red}{\rule[0.025in]{0.05in}{1.5pt}\,\rule[0.025in]{0.05in}{1.5pt}}), and the error function
		$\cbfH(s,p)- \widehat{\cbfH}(s,p)$ (\textcolor{green}{\rule[0.025in]{0.12in}{1.5pt}}) at $p=1.57$. }
		\label{fig:pvf:synthPolApprox}
	\end{figure}

\vspace{-6ex}
	\begin{figure}[H]	
		\centering\begin{minipage}{0.49\textwidth}
			\ifuseTikzGraphs
				\begin{adjustbox}{max width=0.9\textwidth}
				\pgfplotsset{ymin=0.001, ymax=10, xlabel={Frequency $\omega$},ylabel={Magnitude},title={}}
				\input{"Graphics/PVF_SynthP/PolyApprox_r4_sPt6.tex"}
			\end{adjustbox}
			\else\includegraphics[width=0.9\textwidth, trim= 0 0 0 0, clip]{PVF_SynthP/PolyApprox_r4_sPt6}\fi
			\subcaption{$p=3.86$ -- \Cref{alg:pvf:paramVF1}}
		\end{minipage} \hfill
		\begin{minipage}{0.49\textwidth}
			\ifuseTikzGraphs
				\begin{adjustbox}{max width=0.9\textwidth}
				\pgfplotsset{ymin=0.00031, ymax=10, xlabel={Frequency $\omega$},ylabel={Magnitude}}
				\input{"Graphics/PVF_SynthP/RatApprox_r4_sPt6.tex"}
			\end{adjustbox}
			\else\includegraphics[width=0.9\textwidth, trim= 0 0 0 0, clip]{PVF_SynthP/RatApprox_r4_sPt6}\fi
			\subcaption{$p=3.86$ -- \Cref{alg:pvf:paramVFVarProj}}
		\end{minipage}
\caption{Frequency comparisons of $\cbfH(s,p)$ (\textcolor{blue}{\rule[0.025in]{0.12in}{1.5pt}}), $\widehat{\cbfH}(s,p)$ (\textcolor{red}{\rule[0.025in]{0.05in}{1.5pt}\,\rule[0.025in]{0.05in}{1.5pt}}), and the error function
		$\cbfH(s,p)- \widehat{\cbfH}(s,p)$ (\textcolor{green}{\rule[0.025in]{0.12in}{1.5pt}}) at $p=3.86$}
		\label{fig:pvf:synthRatApprox}
	\end{figure}

		To illustrate this more clearly,  we show the 
$\mathcal{H}_2$ and $\mathcal{H}_\infty$
approximation errors over the entire parameter domain $\calP = [1,5]$ in \cref{fig:pvf:synthRatPolComparison}. 
We observe that, except for very few points, adaptive rational basis functions using  \Cref{alg:pvf:paramVFVarProj} yields superior approximation of $\widehat\cbfH(s,p)$.
	\begin{figure}[H]
		\centering
		\ifuseTikzGraphs
		\begin{adjustbox}{width=0.85\textwidth}
\begin{tikzpicture}

\begin{axis}[%
scale only axis,
xmin=1,
xmax=5,
xlabel style={font=\color{white!15!black}},
xlabel={Parameter},
ymode=log,
ymin=0.01,
ymax=1,
yminorticks=true,
ylabel style={font=\color{white!15!black}},
ylabel={Relative $\cbfH_2$ error},
axis background/.style={fill=none},
title style={font=\bfseries},
title={Approximation Quality over Parameter Range}, %Approximation Quality with 4 basis functions},
unit vector ratio*=1 0.75 1,
xmajorgrids,
ymajorgrids,
yminorgrids,
legend pos = outer north east,
legend style={legend cell align=left, align=left, draw=white!15!black}
]
\addplot [color=blue,line width=2pt]
  table[row sep=crcr]{%
1	0.0530283506824558\\
1.04040404040404	0.0433081033503691\\
1.08080808080808	0.0362924359933657\\
1.12121212121212	0.0323181755059271\\
1.16161616161616	0.0313168249729593\\
1.2020202020202	0.0326069330190318\\
1.24242424242424	0.0379901217729351\\
1.28282828282828	0.0438432864390484\\
1.32323232323232	0.0490337882941489\\
1.36363636363636	0.0534893570198831\\
1.4040404040404	0.0571989359507645\\
1.44444444444444	0.0601831750918594\\
1.48484848484848	0.06248264310716\\
1.52525252525253	0.0641534035755116\\
1.56565656565657	0.0652660354711901\\
1.60606060606061	0.0659064043006442\\
1.64646464646465	0.0695133597247378\\
1.68686868686869	0.0749815615630354\\
1.72727272727273	0.0813349680063237\\
1.76767676767677	0.088607962179642\\
1.80808080808081	0.0967847756877455\\
1.84848484848485	0.105804130956099\\
1.88888888888889	0.115565267217299\\
1.92929292929293	0.125933491534957\\
1.96969696969697	0.136744430922416\\
2.01010101010101	0.147806967924939\\
2.05050505050505	0.158905279667093\\
2.09090909090909	0.169800554100057\\
2.13131313131313	0.180232957962123\\
2.17171717171717	0.189924373760142\\
2.21212121212121	0.198582356848129\\
2.25252525252525	0.205905707385216\\
2.29292929292929	0.211592012176518\\
2.33333333333333	0.215347495645891\\
2.37373737373737	0.216899547122039\\
2.41414141414141	0.216012409036457\\
2.45454545454545	0.212506812978429\\
2.49494949494949	0.206285026439847\\
2.53535353535354	0.197364182964197\\
2.57575757575758	0.185923563469956\\
2.61616161616162	0.172376661513468\\
2.65656565656566	0.157486859943186\\
2.6969696969697	0.142550514612739\\
2.73737373737374	0.129636643449472\\
2.77777777777778	0.129844155661514\\
2.81818181818182	0.131528515352883\\
2.85858585858586	0.134639454449875\\
2.8989898989899	0.150917582468266\\
2.93939393939394	0.176811727288776\\
2.97979797979798	0.207214490682203\\
3.02020202020202	0.240323318867847\\
3.06060606060606	0.274826102369814\\
3.1010101010101	0.309757675862035\\
3.14141414141414	0.344382640544936\\
3.18181818181818	0.378118726964137\\
3.22222222222222	0.410487205473071\\
3.26262626262626	0.44107794974602\\
3.3030303030303	0.469521177139067\\
3.34343434343434	0.495460903104466\\
3.38383838383838	0.518526669559209\\
3.42424242424242	0.538300675087987\\
3.46464646464646	0.554277551045079\\
3.50505050505051	0.565814272900068\\
3.54545454545455	0.572069111574001\\
3.58585858585859	0.571933316666764\\
3.62626262626263	0.563971859089618\\
3.66666666666667	0.539226139717426\\
3.70707070707071	0.485831581368534\\
3.74747474747475	0.426674705370437\\
3.78787878787879	0.363462018249714\\
3.82828282828283	0.29914470764668\\
3.86868686868687	0.239211698805744\\
3.90909090909091	0.194282156749187\\
3.94949494949495	0.180607526813533\\
3.98989898989899	0.205712598513095\\
4.03030303030303	0.257185497182535\\
4.07070707070707	0.301806745366282\\
4.11111111111111	0.34030286963964\\
4.15151515151515	0.371475164238472\\
4.19191919191919	0.394990612326425\\
4.23232323232323	0.411370364004384\\
4.27272727272727	0.42135382524861\\
4.31313131313131	0.425652235728024\\
4.35353535353535	0.42486493905057\\
4.39393939393939	0.419461469816585\\
4.43434343434343	0.409788031723163\\
4.47474747474747	0.39608094046199\\
4.51515151515152	0.378480050051481\\
4.55555555555556	0.357039708124949\\
4.5959595959596	0.331736780178\\
4.63636363636364	0.30247638061432\\
4.67676767676768	0.269097302121747\\
4.71717171717172	0.231382778694962\\
4.75757575757576	0.189095586041571\\
4.7979797979798	0.142119062481716\\
4.83838383838384	0.0912122378829334\\
4.87878787878788	0.0590520170196573\\
4.91919191919192	0.0633181461657443\\
4.95959595959596	0.133333278858778\\
5	0.220031290982724\\
};
\addlegendentry{Polynomial}

\addplot [color=red,line width=2pt]
  table[row sep=crcr]{%
1	0.0328697175763129\\
1.04040404040404	0.0263504214576962\\
1.08080808080808	0.0204676256464274\\
1.12121212121212	0.02321476128297\\
1.16161616161616	0.0303675264910665\\
1.2020202020202	0.040029623952681\\
1.24242424242424	0.0493765449614462\\
1.28282828282828	0.0576838714537002\\
1.32323232323232	0.0649075902154857\\
1.36363636363636	0.0710189962352715\\
1.4040404040404	0.0759984006836552\\
1.44444444444444	0.0798332557971924\\
1.48484848484848	0.0825179974527079\\
1.52525252525253	0.0840548920587584\\
1.56565656565657	0.0844557809936158\\
1.60606060606061	0.0837449057521793\\
1.64646464646465	0.0819632553097034\\
1.68686868686869	0.0791752359730742\\
1.72727272727273	0.075479043839034\\
1.76767676767677	0.0710230533343894\\
1.80808080808081	0.0660318320911978\\
1.84848484848485	0.0608462309656701\\
1.88888888888889	0.0559781666332095\\
1.92929292929293	0.0521588522116462\\
1.96969696969697	0.0503003990651348\\
2.01010101010101	0.0512452012612225\\
2.05050505050505	0.0553614113972523\\
2.09090909090909	0.0623803505191975\\
2.13131313131313	0.0716486334856469\\
2.17171717171717	0.0824531855406119\\
2.21212121212121	0.0941697515244999\\
2.25252525252525	0.10627743790867\\
2.29292929292929	0.118331736893446\\
2.33333333333333	0.129937933761095\\
2.37373737373737	0.140734339395324\\
2.41414141414141	0.150384562447337\\
2.45454545454545	0.15857632593526\\
2.49494949494949	0.165024678180762\\
2.53535353535354	0.169477974935885\\
2.57575757575758	0.171725363076945\\
2.61616161616162	0.171604699593554\\
2.65656565656566	0.169009980943226\\
2.6969696969697	0.163897521836564\\
2.73737373737374	0.156290384803317\\
2.77777777777778	0.146281006556012\\
2.81818181818182	0.134032741499624\\
2.85858585858586	0.119782523335787\\
2.8989898989899	0.103850176121593\\
2.93939393939394	0.0950323302298873\\
2.97979797979798	0.0890339452876156\\
3.02020202020202	0.0824152236375095\\
3.06060606060606	0.0753622731152703\\
3.1010101010101	0.0680960398027433\\
3.14141414141414	0.0622680146521131\\
3.18181818181818	0.0722642011252617\\
3.22222222222222	0.0816602374445337\\
3.26262626262626	0.09363426407999\\
3.3030303030303	0.111739935273005\\
3.34343434343434	0.128832201732473\\
3.38383838383838	0.144514832766811\\
3.42424242424242	0.158439278443366\\
3.46464646464646	0.170260603514212\\
3.50505050505051	0.17960793065741\\
3.54545454545455	0.186059629476778\\
3.58585858585859	0.189121899494945\\
3.62626262626263	0.188216518349174\\
3.66666666666667	0.180289879158068\\
3.70707070707071	0.161445883529153\\
3.74747474747475	0.13959297613674\\
3.78787878787879	0.118021125161927\\
3.82828282828283	0.096974675181649\\
3.86868686868687	0.0818509297875568\\
3.90909090909091	0.0802015708161591\\
3.94949494949495	0.0797897805747511\\
3.98989898989899	0.0811811605810471\\
4.03030303030303	0.0845747261019642\\
4.07070707070707	0.0846410831821342\\
4.11111111111111	0.0872048841538039\\
4.15151515151515	0.091512880139009\\
4.19191919191919	0.09353516409396\\
4.23232323232323	0.0977458327767003\\
4.27272727272727	0.0993862904181674\\
4.31313131313131	0.0982372640015704\\
4.35353535353535	0.0947407316548239\\
4.39393939393939	0.0893267798834366\\
4.43434343434343	0.0824142658531616\\
4.47474747474747	0.0744115948465074\\
4.51515151515152	0.0657206856819276\\
4.55555555555556	0.0567505820852786\\
4.5959595959596	0.0479530660644497\\
4.63636363636364	0.0399026413614763\\
4.67676767676768	0.0334419232125337\\
4.71717171717172	0.0297905153331266\\
4.75757575757576	0.0314167998242569\\
4.7979797979798	0.0427060920989925\\
4.83838383838384	0.0549222418065164\\
4.87878787878788	0.0671641654339266\\
4.91919191919192	0.079447450449748\\
4.95959595959596	0.0918148719285522\\
5	0.104338968379047\\
};
\addlegendentry{Rational}

\end{axis}
\end{tikzpicture}%
		\end{adjustbox}
		\else\includegraphics[width=0.8\textwidth]{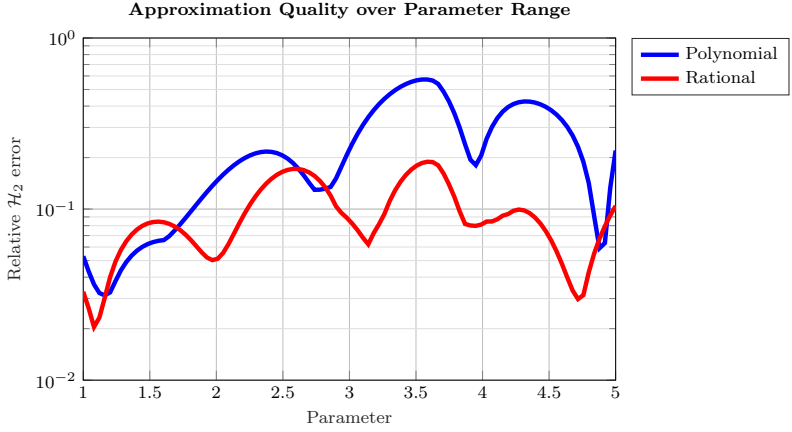}\fi
		\caption{Relative $\mathcal{H}_2$ and $\mathcal{H}_\infty$  errors for \cref{alg:pvf:paramVFVarProj} and 
		\cref{alg:pvf:paramVF1}
		}
		\label{fig:pvf:synthRatPolComparison}
	\end{figure}

\end{example}

\begin{example}\label{ex:pvf:beamPolRat1}
	We  revisit the beam model from \cref{ex:pvf:polBeamExample}. We used the same samples in $s$ and $p$. 
However, in contrast to \cref{ex:pvf:polBeamExample}, here we apply \Cref{alg:pvf:paramVFVarProj} using rational basis functions with  order $r_p = 6$.

	\begin{figure}[H]
		\begin{minipage}{0.49\textwidth}
			\ifuseTikzGraphs
			\begin{adjustbox}{max width=\textwidth}
				\pgfplotsset{ymin=1e-6, ymax=1e+2, xlabel={Frequency $\omega$},ylabel={Magnitude}}
				\input{Graphics/PVF_1P/BeamRev_rs10_Rat_NS1.tex} 
			\end{adjustbox}
			\else\includegraphics[width=0.9\textwidth, trim= 0 0 0 25, clip]{PVF_1P/BeamRev_rs10_Rat_NS1}\fi
			\subcaption{$p=0.01$}
		\end{minipage}\hfill\begin{minipage}{0.49\textwidth}
			\ifuseTikzGraphs
			\begin{adjustbox}{max width=\textwidth}
				\pgfplotsset{ymin=1e-6, ymax=1e+2, xlabel={Frequency $\omega$},ylabel={Magnitude}}
				\input{Graphics/PVF_1P/BeamRev_rs10_Rat_NS2.tex}
			\end{adjustbox}
			\else\includegraphics[width=0.9\textwidth, trim= 0 0 0 25, clip]{PVF_1P/BeamRev_rs10_Rat_NS2}\fi
			\subcaption{$p=0.02$}
		\end{minipage}

	\caption{Frequency responses of $\cbfH(s,p)$ (\textcolor{blue}{\rule[0.025in]{0.12in}{1.5pt}}), $\widehat{\cbfH}(s,p)$ (\textcolor{red}{\rule[0.025in]{0.05in}{1.5pt}\,\rule[0.025in]{0.05in}{1.5pt}}), and the error function
		$\cbfH(s,p)- \widehat{\cbfH}(s,p)$ (\textcolor{green}{\rule[0.025in]{0.12in}{1.5pt}}) at $p=0.01$ and $p=0.22$ using rational basis functions}
		\label{fig:pvf:errorPlotSampledPointsBeam}
	\end{figure}
	The amplitude plots of $\cbfH(s,p)$ and $\widehat{\cbfH}(s,p)$ due to \cref{alg:pvf:paramVFVarProj} are
	shown, in \cref{fig:pvf:errorPlotSampledPointsBeam}, for two representative sampling points, illustrating the accuracy of the approximation. To give a better overall picture, once again we show the approximation error over the continuous parameter interval in \cref{fig:pvf:compRatLocalVFPCont}, 
	illustrating an approximation with a relative error of order $10^{-3}$ for the whole parameter space.

		\begin{figure}[H]
		\centering
		\ifuseTikzGraphs
		\begin{adjustbox}{max width=0.6\textwidth}
\definecolor{mycolor1}{rgb}{0.00000,0.44700,0.74100}%
\definecolor{mycolor2}{rgb}{0.85000,0.32500,0.09800}%
\begin{tikzpicture}

\begin{axis}[%
scale only axis,
xmin=0,
xmax=0.8,
xlabel style={font=\color{white!15!black}},
xlabel={Parameter},
ymode=log,
ymin=3.97480378398524e-05,
ymax=0.00270108279861484,
yminorticks=true,
ylabel style={font=\color{white!15!black}},
ylabel={Error},
axis background/.style={fill=white},
title style={font=\bfseries},
title={Approximation Quality over Parameter Range for Rational $P_\ell$},
xmajorgrids,
ymajorgrids,
yminorgrids,
legend style={legend cell align=left, align=left, draw=white!15!black}
]
\addplot [color=mycolor1, line width=2.0pt]
  table[row sep=crcr]{%
0	0.000701173107839757\\
0.00808080808080808	0.000368713068993989\\
0.0161616161616162	0.000210469937924413\\
0.0242424242424242	0.000161129645663097\\
0.0323232323232323	0.000189981610472579\\
0.0404040404040404	0.000235638536553586\\
0.0484848484848485	0.000274248733627095\\
0.0565656565656566	0.000301842457274769\\
0.0646464646464646	0.000319123023789947\\
0.0727272727272727	0.000327813215457322\\
0.0808080808080808	0.000329764377196728\\
0.0888888888888889	0.00032667215248172\\
0.096969696969697	0.000320035842076483\\
0.105050505050505	0.000311100133101489\\
0.113131313131313	0.000300861688141247\\
0.121212121212121	0.00029005289688054\\
0.129292929292929	0.000279179704692883\\
0.137373737373737	0.000268516476390979\\
0.145454545454545	0.000258151632032837\\
0.153535353535354	0.000248023344927565\\
0.161616161616162	0.000237971799221561\\
0.16969696969697	0.000227786177825163\\
0.177777777777778	0.000217242125302842\\
0.185858585858586	0.000206143636960031\\
0.193939393939394	0.000194337319629487\\
0.202020202020202	0.000181736764404396\\
0.21010101010101	0.000168335208171873\\
0.218181818181818	0.000154197361942196\\
0.226262626262626	0.000139474088007335\\
0.234343434343434	0.00012442428771496\\
0.242424242424242	0.000109433429055291\\
0.250505050505051	9.50581648111918e-05\\
0.258585858585859	8.21127674783757e-05\\
0.266666666666667	7.1741784367985e-05\\
0.274747474747475	6.53519191718613e-05\\
0.282828282828283	6.40540291811988e-05\\
0.290909090909091	6.7805553807018e-05\\
0.298989898989899	7.53883564565227e-05\\
0.307070707070707	8.51974414317363e-05\\
0.315151515151515	9.59764168534958e-05\\
0.323232323232323	0.000106832273639933\\
0.331313131313131	0.000117183221676286\\
0.339393939393939	0.000126657693927231\\
0.347474747474747	0.000135015231259972\\
0.355555555555556	0.000142091672758116\\
0.363636363636364	0.000147824564946145\\
0.371717171717172	0.000152192103205337\\
0.37979797979798	0.000155217969197165\\
0.387878787878788	0.000156975027371737\\
0.395959595959596	0.000157577858008209\\
0.404040404040404	0.000157187592322572\\
0.412121212121212	0.000155992792477784\\
0.42020202020202	0.000154212704284038\\
0.428282828282828	0.000152088982154937\\
0.436363636363636	0.000149862015974777\\
0.444444444444444	0.000147781878986187\\
0.452525252525253	0.000146064081306191\\
0.460606060606061	0.000144880728523589\\
0.468686868686869	0.000144337906870027\\
0.476767676767677	0.000144464201696216\\
0.484848484848485	0.000145199782671738\\
0.492929292929293	0.000146416421124203\\
0.501010101010101	0.00014791245419614\\
0.509090909090909	0.000149439496042726\\
0.517171717171717	0.000150748138674682\\
0.525252525252525	0.000151549148816898\\
0.533333333333333	0.000151595090277903\\
0.541414141414141	0.00015065083317048\\
0.54949494949495	0.000148491237287694\\
0.557575757575758	0.000144957393127715\\
0.565656565656566	0.000139913482655135\\
0.573737373737374	0.000133279223292604\\
0.581818181818182	0.000124992863974821\\
0.58989898989899	0.000115085112480131\\
0.597979797979798	0.00010364142355736\\
0.606060606060606	9.0849207418839e-05\\
0.614141414141414	7.7036742862085e-05\\
0.622222222222222	6.28752541431336e-05\\
0.63030303030303	4.97183481273368e-05\\
0.638383838383838	4.0498064800595e-05\\
0.646464646464647	3.97480378398524e-05\\
0.654545454545455	4.87345866306382e-05\\
0.662626262626263	6.36907944618353e-05\\
0.670707070707071	8.12159327060006e-05\\
0.678787878787879	9.94564542360419e-05\\
0.686868686868687	0.000117372806185669\\
0.694949494949495	0.000134209043465087\\
0.703030303030303	0.000149342850677215\\
0.711111111111111	0.000162186818411203\\
0.719191919191919	0.000172251102799053\\
0.727272727272727	0.000179063076064804\\
0.735353535353535	0.000182295396054895\\
0.743434343434344	0.000181798195246527\\
0.751515151515152	0.000177789094915966\\
0.75959595959596	0.000171218191529023\\
0.767676767676768	0.000164244791522391\\
0.775757575757576	0.000161093748690717\\
0.783838383838384	0.000168164906283083\\
0.791919191919192	0.000191947192589904\\
0.8	0.000235699128145505\\
};
\addlegendentry{$\cbfH_2$-Norm}

\addplot [color=mycolor2, line width=2.0pt]
  table[row sep=crcr]{%
0	0.00270108279861484\\
0.00808080808080808	0.00117847817786706\\
0.0161616161616162	0.000625357742192672\\
0.0242424242424242	0.000433368224020142\\
0.0323232323232323	0.000754623821475184\\
0.0404040404040404	0.000956720366627211\\
0.0484848484848485	0.00107867854606147\\
0.0565656565656566	0.00114595361608784\\
0.0646464646464646	0.00117505718024351\\
0.0727272727272727	0.00117704877424489\\
0.0808080808080808	0.00119398370052834\\
0.0888888888888889	0.00122485836917689\\
0.096969696969697	0.00123519261723243\\
0.105050505050505	0.00122850416225442\\
0.113131313131313	0.00120774820388948\\
0.121212121212121	0.001175424008314\\
0.129292929292929	0.00113364849019022\\
0.137373737373737	0.00108422651137401\\
0.145454545454545	0.00102870005981709\\
0.153535353535354	0.00096838480391264\\
0.161616161616162	0.000904413968629119\\
0.16969696969697	0.000837760597025661\\
0.177777777777778	0.000769260045068553\\
0.185858585858586	0.000699632265844951\\
0.193939393939394	0.000629496921255635\\
0.202020202020202	0.000559391362150588\\
0.21010101010101	0.000489773553215317\\
0.218181818181818	0.000421045516459855\\
0.226262626262626	0.000353551202342509\\
0.234343434343434	0.00028902059252614\\
0.242424242424242	0.000227897776163881\\
0.250505050505051	0.000182597007076289\\
0.258585858585859	0.000151011161410022\\
0.266666666666667	0.000127717685952975\\
0.274747474747475	0.000127008217443446\\
0.282828282828283	0.00015603391701041\\
0.290909090909091	0.000184294607078428\\
0.298989898989899	0.00021344617682968\\
0.307070707070707	0.000243218930141367\\
0.315151515151515	0.000271069899207677\\
0.323232323232323	0.00029696159527102\\
0.331313131313131	0.000333487828500294\\
0.339393939393939	0.000368163149782078\\
0.347474747474747	0.000400170455576224\\
0.355555555555556	0.000429521855840667\\
0.363636363636364	0.00045623419087276\\
0.371717171717172	0.000480332584030599\\
0.37979797979798	0.000501850147723097\\
0.387878787878788	0.00052082105015415\\
0.395959595959596	0.000537284724194314\\
0.404040404040404	0.000551284079355593\\
0.412121212121212	0.000562864764964702\\
0.42020202020202	0.000572073971874572\\
0.428282828282828	0.000578962691344475\\
0.436363636363636	0.000583584539424612\\
0.444444444444444	0.000585990508989909\\
0.452525252525253	0.000586235906816198\\
0.460606060606061	0.000584377312089513\\
0.468686868686869	0.000580471539835993\\
0.476767676767677	0.000574574910240416\\
0.484848484848485	0.000566747027929355\\
0.492929292929293	0.000557045657979673\\
0.501010101010101	0.000545528608135635\\
0.509090909090909	0.000532255951803047\\
0.517171717171717	0.000517287349628988\\
0.525252525252525	0.000500680751632996\\
0.533333333333333	0.00048249579174349\\
0.541414141414141	0.000462791144021265\\
0.54949494949495	0.000441626690322395\\
0.557575757575758	0.000419059555773164\\
0.565656565656566	0.000395149563323577\\
0.573737373737374	0.000369953813355549\\
0.581818181818182	0.000343530571551394\\
0.58989898989899	0.000315937058364598\\
0.597979797979798	0.000287231333820538\\
0.606060606060606	0.00025746921066103\\
0.614141414141414	0.000226708486338584\\
0.622222222222222	0.000195007225091694\\
0.63030303030303	0.000162425386792647\\
0.638383838383838	0.000129027102298389\\
0.646464646464647	9.52818707849855e-05\\
0.654545454545455	6.31782444351149e-05\\
0.662626262626263	6.7047248492777e-05\\
0.670707070707071	8.5149535235992e-05\\
0.678787878787879	0.000102779823034637\\
0.686868686868687	0.000119248993360475\\
0.694949494949495	0.000134080024792261\\
0.703030303030303	0.000164834625805352\\
0.711111111111111	0.000202707710974906\\
0.719191919191919	0.000241125265575385\\
0.727272727272727	0.000280481510189919\\
0.735353535353535	0.000320099242443294\\
0.743434343434344	0.000359933043580097\\
0.751515151515152	0.000399939010594454\\
0.75959595959596	0.000440076054174878\\
0.767676767676768	0.000480302018168542\\
0.775757575757576	0.000520577901259946\\
0.783838383838384	0.000560865471494407\\
0.791919191919192	0.000601125657944036\\
0.8	0.000641323841832755\\
};
\addlegendentry{$\cbfH_\infty$-Norm}

\end{axis}
\end{tikzpicture}%
		\end{adjustbox}
		\else\includegraphics[width=0.8\textwidth,height=0.4\textwidth]{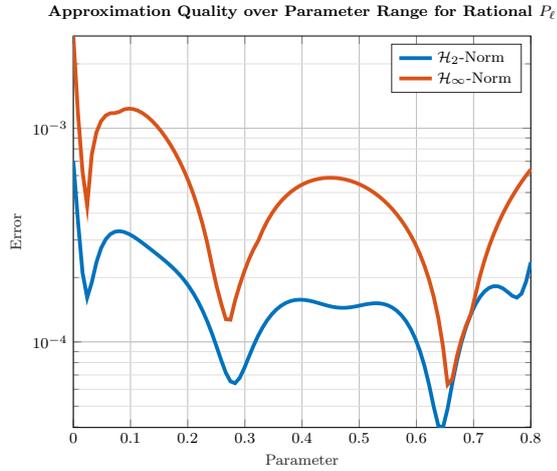}\fi
		\caption{Relative $\cbfH_2$ and $\cbfH_\infty$ errors over $\mathcal{P}$ for rational basis functions with $r_p = 6$.}
		\label{fig:pvf:compRatLocalVFPCont}
	\end{figure}

\end{example}
\begin{example}\label{ex:pvf:beamMeasOutliers}
In this example, we illustrate outlier-resilience of the least squares approximation for fitting parametrized dynamical systems in the framework we have developed.
We used the model and setup as in \cref{ex:pvf:polBeamExample} and \cref{ex:pvf:beamPolRat1}. In order to construct outliers, we chose poor local models at two of the ten sampling points.
	\Cref{fig:PVF:errorCompBeam} shows the performance of \cref{alg:pvf:paramVF1} and \cref{alg:pvf:paramVFVarProj} at the sampling points together with the accuracy of the local \textsf{VF} approximants. Observe that both least-squares soluions
have a nearly uniform error across the parameter domain and avoids the outlier. If one \emph{interpolates} among the local models, instead, the resulting parametric reduced model will suffer large excursions due to the outliers.
\begin{figure}[H]
	\centering
	\ifuseTikzGraphs
	\begin{adjustbox}{max width=0.65\textwidth,center}
\begin{tikzpicture}

\begin{axis}[%
scale only axis,
xmin=0,
xmax=0.8,
xtick={0.1,0.2,...,0.8},
xlabel style={font=\color{white!15!black}},
xlabel={Parameter},
ymode=log,
ymin=1e-06,
ymax=1,
yminorticks=true,
ylabel style={font=\color{white!15!black}},
ylabel={$\mathbf{\mathcal{H}}_2$ Error},
axis background/.style={fill=none},
title style={font=\bfseries},
unit vector ratio*=1 0.2 1,
title={Error Comparison with Measurement Outliers},
xmajorgrids,
ymajorgrids,
yminorgrids,
legend pos = north east,
legend style={legend cell align=left, align=left, draw=white!15!black, style={at={(0.83,0.85)},anchor=west}}
]
\addplot [color=blue, line width=2.0pt, mark=asterisk, mark size = 4pt, mark options={solid, blue}]
  table[row sep=crcr]{%
0.01	0.000101554821002259\\
0.0977777777777778	1.48215752233708e-05\\
0.185555555555556	0.0746274975616341\\
0.273333333333333	3.10066875001578e-06\\
0.361111111111111	0.0657450243194737\\
0.448888888888889	3.35975548090432e-05\\
0.536666666666667	9.87287674783653e-06\\
0.624444444444444	2.92228695578135e-06\\
0.712222222222222	1.68778885346967e-06\\
0.8	0.000125453255202831\\
};
\addlegendentry{Local VF}

\addplot [color=red, line width=2.0pt, mark=asterisk, mark size = 4pt, mark options={solid, red}]
  table[row sep=crcr]{%
0.01	0.0011601430239261\\
0.0977777777777778	0.000761654685208006\\
0.185555555555556	0.000742582445523997\\
0.273333333333333	0.000677534813723023\\
0.361111111111111	0.000677245388519313\\
0.448888888888889	0.0005499403814816\\
0.536666666666667	0.000292020788653047\\
0.624444444444444	0.00042482375736937\\
0.712222222222222	0.000734536330462867\\
0.8	0.00100681502990453\\
};
\addlegendentry{Polynomial $P_\ell$}

\addplot [color=green, line width=2.0pt, mark=asterisk, mark size = 4pt, mark options={solid, green}]
  table[row sep=crcr]{%
0.01	0.00174719355890467\\
0.0977777777777778	0.000688450553269804\\
0.185555555555556	0.000946663868743276\\
0.273333333333333	0.000982580572678576\\
0.361111111111111	0.000806753827777604\\
0.448888888888889	0.000718855891511691\\
0.536666666666667	0.000467575698466038\\
0.624444444444444	0.000314814071529305\\
0.712222222222222	0.000752258900108363\\
0.8	0.00103044928993429\\
};
\addlegendentry{Rational $P_\ell$}

\end{axis}
\end{tikzpicture}%
	\end{adjustbox}
	\else\includegraphics[width=0.65\textwidth]{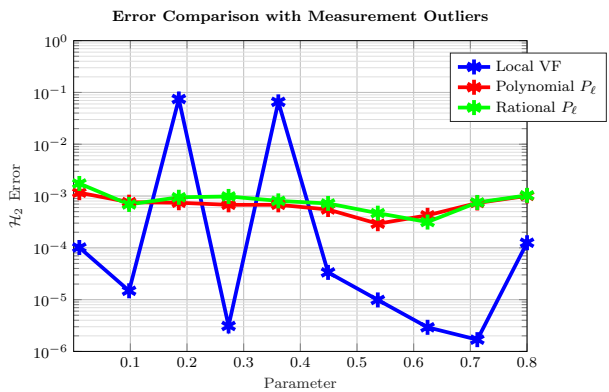}\fi
	\caption{Comparison of relative $\cbfH_2$ errors at sampled parameter values.}
	\label{fig:PVF:errorCompBeam}
\end{figure}

\end{example}

Our observations so far can be summarized as follows:
\begin{enumerate}	
	\item Using rational functions that are adapted to the data for the parametric basis appears to produce approximation quality that is at least as good as what may be obtained by using polynomial functions for the parametric basis, independent of the structure of the original model, see \cref{ex:pvf:beamPolRat1}.
	\item If the underlying model has, in fact, a rational parametrization, then using data-adapted rational functions for the parametric basis may be expected to outperform polynomial bases, see \cref{ex:pvf:synthRational}.
	\item Combining local models with a least-squares measure helps improve resilience with respect to outliers in the data, caused in the example that we considered by a few poor local models.
	\end{enumerate}

%  +  +  +  +  +  +  +  +  +  +  +  +  +  +  +  +  +  +  +  +  +  +  +  +  +  +  +  +  +  +  +  +  +  +  +  +  +  +  +  +  +  +  +  +  +  +  +  +  +
\paragraph{{Real Parametric Systems:}}
For a  rational transfer function, $\cbfH(s)$, to be \emph{real}, it must have a real realization, say, $\cbfH(s) = \bfc^\top(s\bfI-\bfA)^{-1}\bfb$ for $\bfA \in \mathbb{R}^{n\times n}$ and $\bfb,\bfc \in \mathbb{R}^{n}$.  This can be assured if $\cbfH(\overline{s})= \overline{\cbfH(s)}$ for all $s\in \mathbb{C}$ where $\cbfH(s)$ is defined. Analogously we find for the parametrized case that $\cbfH(s,p)$ is \emph{real} if  (a) $\calP$ is closed under conjugation ($\overline{\calP}=\calP$), and 
(b) $\cbfH(\overline{s},\overline{p}) = \overline{\cbfH(s,p)}$ for all $s\in\bbC$, $p\in\calP$.
The local models, $ \localH_k(s)$, used in constructing the intermediate parametrized model $\widehat\cbfH(s,p)$ in \cref{eqn:Hhatform}
 have been generated so that $\overline{\localH_k(s)} = \localH_k(\overline{s})$, for $k=1,\ldots,m_p$.  In order to guarantee that $\widehat\cbfH(s,p)$ is also real, we investigate the parametric dependence. 
 
Let  $\rho:\calI\to\calI$ be a permutation of the index set, $\calI \vcentcolon= \{1,\ldots,r_p\}$, so that $P_\ell(\overline{p}) = \overline{P_{\rho_{\ell}}(p)}$, for $p\in\calP$ and each $\ell=1,\ldots,r_p$.  If $\overline{\hat{x}_{k,\ell}} = \hat{x}_{k,\rho_{\ell}}$ for $k=1,\ldots,m_p$ and $\ell=1,\ldots,r_p$ then a direct calculation shows $\cbfH(s,p)$ must be \emph{real}.  %$$

On the other hand, if $\widehat\cbfH(s,p)$ is real then
$$
\widehat\cbfH(\overline{s},\overline{p})  - \overline{\widehat\cbfH(s,p)}  = \sum_{k=1}^{m_p} \overline{\localH_k(s)} \sum_{\ell=1}^{r_p} ( \hat{x}_{k,\ell} - \overline{\hat{x}_{k,\rho_{\ell}}} )
\overline{P_{\rho_{\ell}}(p)} = 0.
$$ 
Since the  sets $\{\localH_k\}$ and $\{P_\ell(p)\}$ are linearly independent, 
we conclude $\overline{\hat{x}_{k,\ell}} = \hat{x}_{k,\rho_{\ell}}$, $k=1,\ldots,m_p$ and $\ell=1,\ldots,r_p$ and so, 
$\widehat\cbfH(s,p)$  in \cref{eqn:Hhatform} is real if and only if $\overline{\hat{x}_{k,\ell}} = \hat{x}_{k,\rho_{\ell}}$ for $k=1,\ldots,m_p$ and $\ell=1,\ldots,r_p$.

As a consequence, if the parametric basis functions are polynomials with real coefficients, then $\hat{x}_{k,\ell}$ must be real for all $\{k,\ell\}$. If the
 parametric basis functions are elementary rational functions, $P_\ell(s)=\frac{1}{s-\pi_{\ell}}$, then the associated poles, $\{\pi_1,..., \pi_{r_p}\}$ must be closed under conjugation. That is, there is a permutation $\rho$ such that $\overline{\pi_{\ell}} = \pi_{\rho_{\ell}}$, and as a consequence
 $P_\ell(\overline{p}) = \overline{P_{\rho_{\ell}}(p)}$.
As a practical matter, the algorithms that are deployed in solving (\ref{eq:AdaptiveBasisLSsoln1}) or (\ref{eq:AdaptiveBasisLSsoln2}) will not exactly 
preserve the conjugate symmetry, $\overline{\hat{x}_{k,\ell}} = \hat{x}_{k,\rho_{\ell}}$ for all $k,\ell$, and so we enforce the condition explicitly, thus guaranteeing a real system as a final outcome.

%!TEX root = ms.tex
%  +   +  +  +  +  +  +  +  +  +  +  +  +   +  +  +  +  +  +  +  +  +  +  +  +   +  +  +  +  +  +  +  +  +  +  +  +   +  +  +  +  +  +  +  +  +  +  +  +   +  +  +  +  +  +  +  +  +  +  +
\section{Optimal Compression of the Intermediate Parameterized Model} \label{sec:phase2}

In Sections \ref{sec:localmodels} and \ref{sec:pvf:paramVF}, we described how to produce a parameterized intermediate model, 
$\widehat\cbfH(s,p)$, that provides a least squares fit to the given data by combining parametrically localized system models with appropriate parametric bases.  Since the number of local models grows with the number of parameter samples, the $s$-order of $\widehat\cbfH(s,p)$ also grows rapidly  as the number of parameter samples grows.  For example, suppose $m_p = 20$ parameter samples are used with uniform local model order,  $\nu = \nu_k = 50$.  The resulting $s$-order of $\widehat\cbfH(s,p)$ may be as large as $n_s = m_p \cdot \nu = 1000$.  Could we do better ?  

 One might expect that adjacent parameter values could produce system responses that have many common features, leading one to suspect that there could be significant redundancy among the local models. Thus, we might reasonably expect that a lower order system could be found that might also fit the data well.   In this section, we describe a second phase to our solution process that \emph{compresses} $\widehat\cbfH(s,p)$.  We produce another parameterized model, $\widehat\cbfH_{red}(s,p)$, that approximates 
 $\widehat\cbfH(s,p)$ (and hence the data) quite well (indeed, optimally in a sense we specify); it has a parameterization that conforms with the original model parameterization in the sense that for each $\hat{s}\in \mathbb{C}$, $\widehat\cbfH_{red}(\hat{s},\cdot)\in\mathfrak{P}_{r_p}$ as well; but we anticipate that $\widehat\cbfH_{red}(s,p)$ will have $s$-order substantially smaller than $n_s$.
 
Our approach uses an $\mathbf{\mathcal{H}}_2$-optimal model reduction strategy similar to one introduced by Baur et al. \cite{baur2011interpolatory} but specifically tailored for the parametric structure of the intermediate model obtained in the first phase of our solution process described in \Cref{sec:pvf:paramVF}.   The compressed model coming out of this second phase retains the parametric structure of the intermediate model while producing a locally optimal approximation with respect to a continuous system-theoretic error measure which we describe below. 

\subsection{Optimal compression with respect to an ${\calH_2\otimes L_2(\calP)}$ error}
For each parameter sample $\mu_k$, $k=1,\ldots,m_p$, the corresponding local model,
$\localH_k(s)$, can be realized as
\begin{equation}
	 \localH_k(s) =\bfc_k^\top\left( s \bfI -\bfA_k \right)^{-1}\bfb_k,\qquad k=1,\ldots,m_p,
\end{equation}
where $\bfA_k \in \bbC^{\nu_k \times \nu_k}$, and $\bfb_k,\bfc_k \in \bbC^{\nu_k}$.  
Define,
\begin{equation}\label{eq:pvf:defSysMatricesPVFIRKA}
	\mathfrak{A} \vcentcolon= \mathsf{diag}(\bfA_1, \bfA_2,..., \bfA_{m_p}) =
	\begin{bmatrix}
		\bfA_1 \\ & \ddots \\ &&\bfA_{m_p}
	\end{bmatrix}\in\bbC^{n_s\times n_s},\qquad
	\mathfrak{B} \vcentcolon=
	\begin{bmatrix}
		\widehat\bfb_1 \\ \vdots \\ \widehat\bfb_{m_p}
	\end{bmatrix}\in\bbC^{n_s}.
\end{equation}
and the block diagonal matrix:
$$\mathbb{D}=\mathsf{diag}(\bfc_1, \bfc_2,..., \bfc_{m_p})
=\left[\begin{array}{cccc} \bfc_1 & 0 & \ldots & 0 \\
0 & \bfc_2  &   &  \\
\vdots &  & \ddots &  \\
0 &  & & \bfc_{m_p} 
\end{array}\right] \in \bbC^{n_s \times m_p}$$
(note that $\mathbb{D}$ is not a square matrix and moreover each $\bfc_k$ will be of different order generally).
For the parametric basis defining the factor space, $\mathfrak{P}_{r_p}$ (or $\mathfrak{P}_{r_p}^\bfpi$ if it was chosen adaptively as in \Cref{sec:pvf:paramVF}),  define
\begin{equation}\label{eq:pvf:defBBV}
	\bbV(p) \vcentcolon = \begin{bmatrix}P_1(p) & \ldots & P_{r_p}(p) \end{bmatrix}^\top \in \bbC^{r_p}
\end{equation}
(We suppress the meta-parameter $\hat{\bfpi}$ if a data-adapted parameter basis is used).
Then, the parametric model resulting from either \cref{alg:pvf:paramVF1} or \ref{alg:pvf:paramVFVarProj} may be written as
\begin{equation}\label{eq:pvf:combLocalVF}
\widehat\cbfH(s,p)=\sum_{k,\ell} \hat{x}_{k\ell} \,\mathfrak{h}_k(s)P_{\ell}(p) 
= \sum_{k,\ell} P_{\ell}(p) \hat{x}_{k\ell} \,\bfc_k^\top\left( s \bfI -\bfA_k \right)^{-1}\bfb_k  
= \bbV(p)^\top  \mathbb{G}(s),
\end{equation} 
where $ \mathbb{G}(s)=(\mathbb{D}\widehat{X})^\top(s\bfI-\mathfrak{A})^{-1}\mathfrak{B}$.   Notice that for any $s\in \mathbb{C}$, $\mathbb{G}(s)\in \bbC^{r_p}$, which can then be viewed as a \emph{parameter-free} SIMO (single-input/multiple-output) system mapping scalar inputs to $r_p$-dimensional vector outputs.  

We seek a parameterized model, $\widehat\cbfH_{red}(s,p)$, with a compatible parameterization to $\widehat\cbfH(s,p)$ and with (potentially) significantly smaller $s$-order, say, $n_{red}\ll n_s$.   This invites the ansatz that $\widehat\cbfH_{red}(s,p)$ must have the form
$$
\widehat\cbfH_{red}(s,p) = \bbV(p)^\top \mathbb{G}_{red}(s) \qquad \mbox{with}\quad   \mathbb{G}_{red}(s)= C_{red}(s\bfI-A_{red})^{-1}B_{red},
$$
and
\begin{equation}\label{eq:pvf:defSysMatricesRedModel}
	\bfA_{red}\in\bbC^{n_{red}\times n_{red}},\qquad
	\bfB_{red}\in \bbC^{n_{red}},\qquad\mbox{and}\quad \bfC_{red}\in \bbC^{r_p \times n_{red}}.
\end{equation}

The error with which $\widehat\cbfH_{red}(s,p)$ approximates $\widehat\cbfH(s,p)$ can be measured through 
a joint norm defined for bivariate functions $\cbfH \in \calH_2\otimes L_2(\calP) $ as
\begin{equation}\label{eq:pvf:combSysNormPVFIRKA}
	\| \cbfH \|_{\calH_2\otimes L_2(\calP)} \vcentcolon= \sqrt{\frac{1}{2\pi} \int_{-\infty}^{+\infty}
\iint\limits_{\mathcal{P}} | \cbfH(\rmi \omega,\, p) |^2
\,\rmd A(p)\,\rmd\omega}.
\end{equation} 
where  $\rmd A(\bfparp)$ is either planar Lebesque measure defined on our parameter set, $\calP\subset \mathbb{C}$, or linear Lebesque measure, if $\calP$ is a line segment. 

In particular, note that 
\begin{equation*}
\begin{aligned}
\| \widehat\cbfH - \widehat\cbfH_{red}\|_{\calH_2\otimes L_2(\calP)}^2 & =\frac{1}{2\pi} \int_{-\infty}^{+\infty}
\iint\limits_{\mathcal{P}} |\widehat\cbfH(\imath \omega,\, \bfparp)-\widehat\cbfH_{red}(\imath \omega,\, p)|^2\,\rmd A(p)\,\rmd\omega \\
& =\frac{1}{2\pi} \int_{-\infty}^{+\infty}
\iint\limits_{\mathcal{P}} |\bbV(p)^\top\left(\mathbb{G}(\imath \omega)-\mathbb{G}_{red}(\imath \omega)\right)|^2\,\rmd A(p)\,\rmd\omega \\
& =\frac{1}{2\pi} \int_{-\infty}^{+\infty}
\|\mathbf{R}\left(\mathbb{G}(\imath \omega)-\mathbb{G}_{red}(\imath \omega)\right)\|_2^2\,\rmd\omega \\
&= \| \mathbf{R}\mathbb{G}-\mathbf{R}\mathbb{G}_{red}\|_{\calH_2}^2
\end{aligned}
\end{equation*}
where $\mathbf{R}\in \bbC^{r_p \times r_p} $ is the (upper triangular) Cholesky factor of the Gram matrix of the parametric basis with respect to $L_2(\calP)$: 
$\mathbf{R}^\top \mathbf{R}=\iint\limits_{\mathcal{P}} \bbV(p)\bbV(p)^\top\, \rmd A(p)$.  The last line recognizes the usual definition of the $\calH_2$ norm of a SIMO dynamical system. 

This shows that finding a compressed parameterized model, $\widehat\cbfH_{red}(s,p)$, that is an \emph{optimal} approximation to $\widehat\cbfH(s,p)$ with respect to the error measure (\ref{eq:pvf:combSysNormPVFIRKA}) is \emph{equivalent} to finding an
$\calH_2$-optimal reduced order model, $G_{red}(s)$, approximating a (preweighted) system, $\mathbf{R}\mathbb{G}(s)$, and then unweighting the result by defining  $\mathbb{G}_{red}(s)=\mathbf{R}^{-1}G_{red}(s)$.   Notice that once we have determined $\mathbb{G}_{red}(s)$, our optimally compressed parameterized model is available as $\widehat\cbfH_{red}(s,p)= \bbV(p)^\top \mathbb{G}_{red}(s)$. 

The equivalence between an ${\calH_2\otimes L_2(\calP)}$-optimal parametric model reduction problem and a weighted  parameter-free $\calH_2$-optimal model reduction problem is an innovation introduced in \cite{baur2011interpolatory} that allows one to solve ${\calH_2\otimes L_2(\calP)}$ approximation problems by using well-established, numerically efficient tools 
for $\calH_2$ model reduction, e.g., the iterative rational Krylov method (IRKA) of \cite{GAB08}.  While \cite{baur2011interpolatory} considered only affine parameter dependence, we allow more general parametric families.  For more details on the $\calH_2$ optimal model reduction problem and in particular on IRKA, see \cite{GAB08,ABG10}.

We summarize this discussion and describe an overall two-phase procedure below:
\begin{algorithm}[H]
	\caption{Aggregate two-phase algorithm}\label{alg:overall}
	\begin{tabularx}{\textwidth}{lX}
		{\bf INPUT:}&  Measurements $\left\{ \xi_i,\mu_j,\cbfH(\xi_i,\mu_j)\right\}_{i=1,j=1}^{i=m_s,j=m_p}$, \\
		& Target model order: $n_{red}$ \\[1mm]
	{\bf OUTPUT: } & Parametrized model $\widehat\cbfH_{red}(s,p)$ with $s$-order $n_{red}$ 
	\end{tabularx}
	\begin{enumerate}
		\item \textbf{Phase 1}: Apply \Cref{alg:pvf:paramVF1} or  \ref{alg:pvf:paramVFVarProj} to construct an intermediate model
$\widehat\cbfH(s,p)$.
\item \textbf{Phase 2}: $\calH_2$-compression of intermediate model
\begin{enumerate}
\item Compute the Cholesky factor, $\mathbf{R}$, of the Gram matrix associated with the parametric basis from Phase 1.
 \item Find an $\calH_2$-optimal reduced model, $G_{red}(s)$, of order $n_{red}$ approximating  the (preweighted) SIMO model, $\mathbf{R}\mathbb{G}(s)$, where $\mathbb{G}(s)$ is constructed as in (\ref{eq:pvf:combLocalVF}) from local models derived in Phase 1. 
	\end{enumerate}
\item Return a final parametrized reduced model 
$$
\widehat\cbfH_{red}(s,p) \vcentcolon= \bbV(p)^\top \mathbf{R}^{-1}G_{red}(s) 
$$
	\end{enumerate}
\end{algorithm}

\subsection{Asymptotic stability over the entire parameter domain} 
In most applications, the transfer function $\cbfH(s,p)$ is asymptotically stable for every ${p} \in \mathcal{P}$, i.e.,  for each fixed $\hat{p}$, all the poles of  $\cbfH(s,\hat{p})$ have negative real parts. The structure of the derived model $\widehat{\cbfH}(s,p)$ in \cref{eq:pvf:combLocalVF} provides an assurance of asymptotic stability for every ${p} \in \mathcal{P}$. Indeed, if the local models, $\{\localH_k\}_{1}^{m_p}$, are all asymptotically stable, then our parametrized reduced model $\widehat{\cbfH}(s,p)$ resulting from either \cref{alg:pvf:paramVF1} or \ref{alg:pvf:paramVFVarProj} is guaranteed to be asymptotically stable over the entire parameter domain. 
This is one of the main advantages of knitting together asymptotically stable local reduced models with appropriate parametric basis functions (a point that has also been pointed out in \cite{morBauB09}). Here we use \textsf{VF} to construct the local models. Even though  \textsf{VF} is not guaranteed to produce asymptotically stable systems, in almost all cases it is sufficient to add to \textsf{VF} a ``pole-flipping" step, which involves reflecting intermediate unstable poles that may be encountered back to the left-half plane. In all of our examples we enforced stability on all local models, and so we can guarantee asymptotic stability of $\widehat{\cbfH}(s,p)$  in every case. 
This is in contrast to other parametric data-driven approaches  where a two-variable barycentric-form for $\widehat{\cbfH}(s,p)$ is used either to enforce interpolation  as in the parametric-Loewner formulation of \cite{Lefteriu2011,Ionita2014} or to minimize a least-squares error as in \cite{Grivet-talocia2017} (the same performance criterion considered here).  
Allowing poles to vary with $p$ might allow one to cover a wider range of dynamics with lower order realizations --- but the potential cost is that for these approaches asymptotic stability over the entire parameter domain cannot be guaranteed in general, and for any particular parameter $p=\hat{\mu}$, a stability-correcting post-processing might be required.   With our approach, as we illustrate in various examples, a sufficiently rich local basis, $\mathfrak{R}_{m_p}$, allows us to cover wide-ranging dynamics accurately over the entire parameter domain. 
Another advantage of using the form \cref{eq:pvf:combLocalVF},
as revealed in this section,  is that it allows an optimal parametric-reduction of $\widehat{\cbfH}(s,p)$ to reduce the dimension further for the cases of large number of sampling points. An $\mathcal{H}_2 \otimes L_2(\mathcal{P)}$ optimal reduced model
$\widehat{\cbfH}_{red}(s,p)$ must be, by definition, asymptotically stable. One can either enforce this by introducing a pole-flipping step in IRKA as done in \textsf{VF}, or can use a trust-region variant such as \cite{beattie2009trust}.

\begin{example}
We revisit the beam model from
\Cref{ex:pvf:polBeamExample}. 

We chose $200$ logarithmically spaced frequency samples in  $[10^{-3}\rmi, 10^3\rmi]$ and $m_p = 20$ parameter samples, logarithmically spaced in $[0.0001,1]$. In Step 1. of \Cref{alg:overall}, we use 
\Cref{alg:pvf:paramVF1} where  the order of the local models is $\nu_s=25$ and the polynomial order is  $r_p = 20$.
This results in the intermediate parametric model $\widehat\cbfH(s,p)$ of order $n_s = 500$.
Then, in Step 3. of \Cref{alg:overall}, we perform $\mathcal{H}_2 \otimes L_2(\mathcal{P})$ reduction with $n_{red} = 20$.
In \cref{fig:PVF:beamVFIRKAEx1}, we show the frequency response comparison for two representative parameter values of $p=0.001$ and $p=0.22$. The figure illustrates that the reduced model  $\widehat\cbfH_{red}(s,p)$ of \Cref{alg:overall} matches the approximation quality of the reduced model $\widehat\cbfH_F(s,p)$  of
\Cref{alg:pvf:paramVF1} despite having a much smaller  order.

\begin{figure}[H]
	\centering
	\begin{minipage}{0.49\textwidth}
	\ifuseTikzGraphs
		\begin{adjustbox}{max width=\textwidth,center}
		\pgfplotsset{ymin=1e-8, ymax=1e+2, xlabel={Frequency $\omega$},ylabel={Magnitude}}
		\input{"Graphics/PVF_IRKA/Beam1D_rs25_BODECompP1.tex"} 
		\end{adjustbox}
	\else\includegraphics[width=\textwidth]{"Graphics/PVF_IRKA/PVFIRKA_Beam_r20BODEPlot_RedModelPar1"}\fi
		\subcaption{$p=0.001$}
	\end{minipage}
	\begin{minipage}{0.49\textwidth}
	\ifuseTikzGraphs
		\begin{adjustbox}{max width=\textwidth,center}
		\pgfplotsset{ymin=1e-8, ymax=1e+2, xlabel={Frequency $\omega$},ylabel={Magnitude}}
		\input{"Graphics/PVF_IRKA/Beam1D_rs25_BODECompP2.tex"}
		\end{adjustbox}
	\else\includegraphics[width=\textwidth]{"Graphics/PVF_IRKA/PVFIRKA_Beam_r20BODEPlot_RedModelPar5"}\fi
		\subcaption{$p=0.22$}
	\end{minipage}
\caption{Frequency responses of $\cbfH(s,p)$ (\textcolor{blue}{\rule[0.025in]{0.12in}{1.5pt}}), $\widehat{\cbfH}(s,p)$ (\textcolor{red}{\rule[0.025in]{0.05in}{1.5pt}\,\rule[0.025in]{0.05in}{1.5pt}}), and the error function
		$\cbfH(s,p)- \widehat{\cbfH}(s,p)$ (\textcolor{green}{\rule[0.025in]{0.12in}{1.5pt}}) at $p=0.001$ and $p=0.22$}
	\label{fig:PVF:beamVFIRKAEx1}
\end{figure}

To further illustrate the success of the $\calH_2$-compression step in \Cref{alg:overall}, we pick three levels of reduced orders, namely $n_{red} = 10$, $n_{red} = 20$, and  $n_{red} = 30$. For each of these three cases and for the intermediate model, we compute the relative $\mathcal{H}_2$ error for the entire parameter domain. The results, depicted in   \Cref{fig:pvf:beamCompPVFIRKAH2}, show that  with $n_{red} = 20$, and  $n_{red} = 30$, \Cref{alg:overall} matches the approximation quality of \Cref{alg:pvf:paramVF1} over the entire domain. In \cref{fig:pvf:beamCompStagesErr}, for $p=0.001$, we show the frequency response error plots. Observe that as $n_{red}$ increases to $30$, the error due to  \Cref{alg:overall} matches very closely the error due to \Cref{alg:pvf:paramVF1}.

\begin{figure}[H]
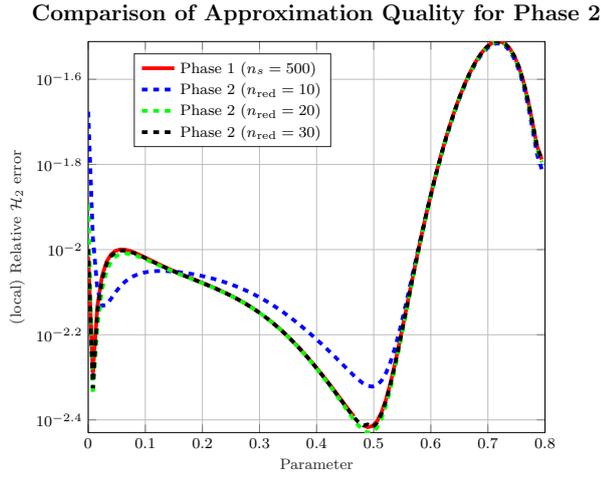

	\centering
	\ifuseTikzGraphs
		\begin{adjustbox}{width=0.65\textwidth,center}
		\input{"Graphics/PVF_IRKA/BeamH2ErrIntervalCompNew.tex"}
		\end{adjustbox}
	\else\includegraphics[width=0.8\textwidth]{PVF_IRKA/BeamH2ErrIntervalCompNew}\fi
	\caption{Relative $\cbfH_2$ error over the entire parameter range $[0,0.8]$. }
	\label{fig:pvf:beamCompPVFIRKAH2}
\end{figure}
\begin{figure}[H]
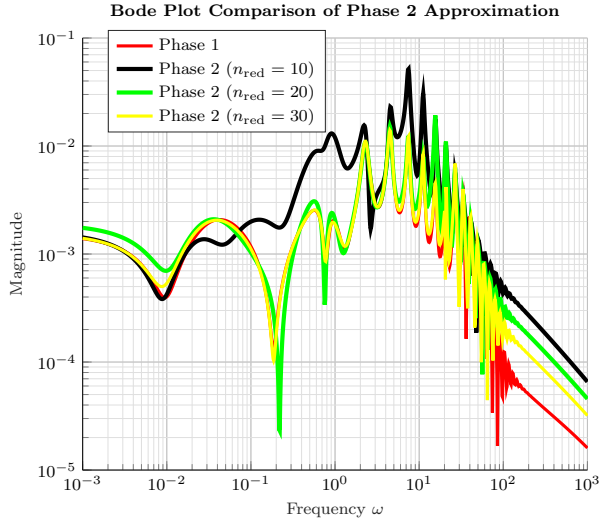

	\centering
	\ifuseTikzGraphs
		\begin{adjustbox}{max width=0.65\textwidth,center}
		\input{"Graphics/PVF_IRKA/Beam1D_rs25_BODEErrCompStagesP1.tex"}
		\end{adjustbox}
	\else\includegraphics[width=0.8\textwidth]{PVF_IRKA/Beam1D_rs25_BODEErrCompStagesP1}\fi
	\caption{Frequency response error for various $n_{red}$ values at $p=0.001$}
	\label{fig:pvf:beamCompStagesErr}
\end{figure}

\end{example}

%!TEX root = ms.tex
\section{Parametric Fitting with Several Parameters}\label{sec:pvf:severalParam}
Up to this point, we have assumed access to frequency response data associated with a single scalar parameterization.  
A moment reflection suggests that the approach we have proposed should have a trivial extension to the case of several parameters,
 $\bfp = \begin{bmatrix} p_1&p_2&\ldots&p_{d}\end{bmatrix}$;
we describe briefly one possible way this extension could proceed for the case of two parameters ($d=2$), labeled as $\bfp = \begin{bmatrix} p & q\end{bmatrix}$. 
We choose parametric basis functions in $p$ and $q$ independently:
Take $P_{k}(p)$, $k=1,\ldots,r_p$ to be basis functions associated with the first parameter $p$ spanning a subspace $\mathfrak{P}_{r_p}$, and $Q_{\ell}(q)$, $\ell=1,\ldots,r_q$ be basis functions associated  with the second parameter $q$ spanning a subspace $\mathfrak{Q}_{r_q}$.
Suppose frequency response observations are taken at frequencies, $s=\xi_1,\ldots, \xi_{m_s}$ for parameter samplings
$p=\mu_1,\ldots, \mu_{m_p}$ and $q=\eta_1,\ldots, \eta_{m_q}$; these observations may be indexed as 
$\mathcal{H}(\xi_i,\mu_{j_1},\eta_{j_2})=\mathbb{H}_{i,(j_1,j_2)}$ with  $i=1,\ldots ,m_s$, $j_1=1,\ldots ,m_p$,  and $j_2=1,\ldots ,m_q$.  
We proceed much as we did in (\ref{eqn:Hhatform}), assuming a separable intermediate model: 
\begin{equation} 
\widehat\cbfH(s,p,q) =\sum_{k=1}^{r_s}\sum_{\ell_1=1}^{r_p}\sum_{\ell_2=1}^{r_q}x_{k,(\ell_1,\ell_2)}\, \mathfrak{h}_k(s)P_{\ell_1}(p) Q_{\ell_2}(q)\in \mathfrak{R}_{r_s} \otimes \mathfrak{P}_{r_p} \otimes \mathfrak{Q}_{r_q}
\end{equation}
and seek an $\widehat{X}=[\hat{x}_{k,(\ell_1,\ell_2)}]$ that solves the least squares problem: 
$$
\widehat{X}=\argmin_{{X}}\sum_{i=1}^{m_s} \sum_{j_1=1}^{m_p}\sum_{j_2=1}^{m_q}  \left| \widehat\cbfH(\xi_i,\mu_{j_1},\eta_{j_2}) - 
\cbfH(\xi_i,\mu_{j_1},\eta_{j_2}) \right|^2
$$
While $\mathbb{H}=[\mathcal{H}(\xi_i,\mu_{j_1},\eta_{j_2})]$ may be viewed as a tensor with (tensor) rank $3$ and 
 dimension $m_s \times m_p \times m_q$, it will be useful to flatten $\mathbb{H}$ along indices associated with parameters, so without changing notation, we think of $\mathbb{H}$  as a two-dimensional array, $\mathbb{H}\in \mathbb{C}^{m_s\times (m_p\cdot m_q)}$, with each row 
 associated with observations at a particular frequency, $\xi$, and column entries stored in $q$-major order (that is, the observations $\mathbb{H}_{i,(j_1,j_2)}$ are stored consecutively in row $i$, with $j_1=1,\ldots , m_p$, $j_2=1, \ldots , m_q$, and with the $q$-index, $j_2$, varying most rapidly).   
 Similarly, we flatten $\widehat{X}$ along indices associated with parameters, thinking of  $\widehat{X}$ as a two-dimensional array, $\widehat{X}\in \mathbb{C}^{r_s\times (r_p\cdot r_q)}$ also without a change in notation. 
 Defining \begin{equation}
\begin{array}{c}\mathbb{A}=[\mathfrak{h}_j(\xi_i)]\in \mathbb{C}^{m_s\times r_s},\quad
\mathbb{B}_p=[P_j(\mu_i)]\in \mathbb{C}^{m_p\times r_p}, \\[2mm]
~~\mbox{and}~~ \mathbb{B}_q=[Q_j(\eta_i)]\in \mathbb{C}^{m_p\times r_p},
\end{array}
\end{equation}
we may reformulate our least squares problem concisely as: 
$$
\widehat{X}=\argmin_{{X}}  \left\| \mathbb{A}X (\mathbb{B}_p\otimes\mathbb{B}_q)^\top - \mathbb{H} \right\|_F^2
$$
One may follow the steps previously discussed for the single parameter case in Sections \ref{sec:localmodels}, \ref{sec:pvf:paramVF}, and \ref{sec:phase2} with minor changes now for the multiparameter case.  In principle, the extension is trivial though potentially tedious.  It is worth a cautionary note, however, that the usual computational issues, common to most parametric model reduction approaches, arising from the need to sample a high-dimensional parameter space will occur here as well. 

The straightforward tensor product/grid-sampling strategy outlined here is an elementary extension of the framework put forward in Sections \ref{sec:localmodels} and \ref{sec:pvf:paramVF}, and is expected to work comparably well in the multiparameter case for a small number of parameters.  However, for even a modest number of parameters more subtle strategies may be necessary.   Note first that the total number of local models that are used may increase dramatically as the number of parameters is increased;  for $d$ parameters, $\bfp = \begin{bmatrix} p_1&p_2&\ldots&p_{d}\end{bmatrix}$, with each parameter sampled at $\hat{m}$ values, say, the number of local models generated is $\hat{m}^d$ which grows explosively as $\hat{m}$ increases if $d$ is large, a  common problem for approaches that make use of local models.  
Suppose the parametric dependence with respect to each parameter is represented uniformly with $\hat{r}$ basis vectors.  Then the total number of unknowns to be determined in Phase 1 of \cref{alg:overall}, is $(\hat{m}\hat{r})^{d}$ (number of elements of $\widehat{X}$).  Moreover, if $\hat{m}^d>m_s$ then the resulting least squares problem is rank deficient, and computational strategies must take this into account, potentially at a significant additional cost.
 Therefore, different sampling approaches such as adaptive sparse sampling or greedy sampling will likely be necessary when $d$ is large; see \cite[Section 3.4]{Benner2016} for a brief discussion. Such an approach leaves the framework of \cref{alg:overall} largely unchanged.  A potential alternative that departs somewhat from the framework of \cref{alg:overall} is a multilevel approach that knits together the Phase 1 and Phase 2 steps, allowing smaller subsets of the $\hat{m}^d$ local models to be hierarchically aggregated.   Strategies such as these that could be suitable for modeling systems with large numbers of parameters will not be pursued further here.  We note that parametric modeling problems involving a large number of parameters often produce staggering computational challenges; we anticipate that strategies such as what we offer in \cref{alg:overall} can play an important role in solving such problems, but are not likely to suffice themselves. 

\begin{example}
In this model, taken from  \cite{baur2011interpolatory},
we consider the convection-diffusion model  on unit square $\Omega = [0,1]\times[0,1]$:
\begin{equation}\label{eq:intro:CDPDE}
	\begin{aligned}
		\dfrac{\partial \phi(t;\bfz)}{\partial t} &=  \Delta \phi(t,\bfz) + \bfp \cdot \nabla \phi(t,\bfz) + b(\bfz) u(t) \qquad &\bfz\in\Omega,~t\in(0,\infty), 
	\end{aligned}
\end{equation}
with homogeneous Dirichlet boundary conditions $\phi(t, \bfz) = 0$, for $\bfz \in \partial \Omega$, where $b(\bfz)$ represents the characteristic function of the domain where the forcing function $u(\cdot)$ acts.
The parameter $\bfp = \begin{bmatrix} p & q\end{bmatrix}^\top$ represents convection in both directions. 
Discretizing \cref{eq:intro:CDPDE} with a finite difference scheme yields
\begin{equation}\label{eq:intro:cdModelSSEquation}
		\begin{aligned}
		\cbfH(s,\bfp) =  \bfc^\top \left( s \bfI - (\bfA_0 + p \bfA_1 + q \bfA_2 \right)  \bfb, 
			\end{aligned}
\end{equation}
where $\bfA_0,\bfA_1,\bfA_2 \in \mathbb{R}^{n\times n}$, and $\bfb,\bfc \in \mathbb{R}^n$ with $n=10000$. 
We sample $\cbfH(s,\bfp)$  with a uniform  $6\times 6$ grid in the parameter space $\Omega$ together  with $100$ frequeny points logarithmically spaced in $[10^2,10^6]$ on the imaginary axis for each of the parameter pair. We use polynomial bases and  choose $r_p = r_q = 12$.
The approximation quality of the two-variable parametric approximation is shown in \cref{fig:pvf:cd2PEx01} for a variety of parameter points, illustrating a high-quality parametric approximant.
\begin{figure}[H]
	\centering
	\ifuseTikzGraphs
		\begin{minipage}{0.46\textwidth}
		\begin{adjustbox}{max width = \textwidth, center}
			\pgfplotsset{ymin = 1e-7, ymax = 5e-3, grid style={dotted,gray!10!black}}
			\input{"Graphics/PVF_2P/CD_p301_55_xi100_r12_r12_P1.tex"}
		\end{adjustbox}
	\end{minipage}
	\hfill\begin{minipage}{0.46\textwidth}
		\begin{adjustbox}{max width = \textwidth, center}
			\pgfplotsset{ymin = 1e-7, ymax = 5e-3, grid style={dotted,gray!10!black}}
			\input{"Graphics/PVF_2P/CD_p301_55_xi100_r12_r12_P3.tex"}
		\end{adjustbox}
	\end{minipage} \\ [2ex]
	\begin{minipage}{0.46\textwidth}
		\begin{adjustbox}{max width = \textwidth, center}
			\pgfplotsset{ymin = 1e-7, ymax = 5e-3, grid style={dotted,gray!10!black}}
			\input{"Graphics/PVF_2P/CD_p301_55_xi100_r12_r12_P5.tex"}
		\end{adjustbox}
	\end{minipage}\hfill\begin{minipage}{0.46\textwidth}
		\begin{adjustbox}{max width = \textwidth, center}
			\pgfplotsset{ymin = 1e-7, ymax = 5e-3, grid style={dotted,gray!10!black}}
			\input{"Graphics/PVF_2P/CD_p301_55_xi100_r12_r12_P6.tex"}
		\end{adjustbox}
	\end{minipage} 
	\else\includegraphics[width=\textwidth]{PVF_2P/CD_p301_55_xi100_r12_r12}\fi
	\caption{Frequency responses of $\cbfH(s,\bfp)$ (\textcolor{blue}{\rule[0.025in]{0.12in}{1.5pt}}), $\widehat{\cbfH}(s,\bfp)$ (\textcolor{red}{\rule[0.025in]{0.05in}{1.5pt}\,\rule[0.025in]{0.05in}{1.5pt}}), and the error function
		$\cbfH(s,\bfp)- \widehat{\cbfH}(s,\bfp)$ (\textcolor{green}{\rule[0.025in]{0.12in}{1.5pt}}) for the two-parameter case at selected samples}
	\label{fig:pvf:cd2PEx01}
\end{figure}
\end{example}

%  +   +  +  +  +  +  +  +  +  +  +  +  +   +  +  +  +  +  +  +  +  +  +  +  +   +  +  +  +  +  +  +  +  +  +  +  +   +  +  +  +  +  +  +  +  +  +  +  +   +  +  +  +  +  +  +  +  +  +  +
\section{Conclusions}

 We have presented a two phase approach to construct a parsimonious parametrized model that fits in a least-squares sense frequency response data arising from observations of a parametrized dynamical system.
 Using parametrized basis functions, the first phase of the proposed algorithm combines local models derived from the observed system response across parameter samplings, in order to solve a coupled least-squares data fitting problem taken with respect to both frequency and parameter samples.
 We consider both fixed parameter bases and varying  parametric bases that have been adapted to the given data.
The second-phase of our approach uses $\cbfH_2$-optimal model reduction strategies to eliminate potential redundancy that may exist among the local models obtained in the parametrized intermediate model from the first phase.  Several examples illustrated the performance of our framework.

%
%
%  +  +  +  +  +  +  +  +  +  +  +  +  +  +  +  +  +  +  +  +  +  +  +  +  +  +  +  +  +  +  +  +  +  +  +  +  +  +  +  +  +  +  +  +  +  +  +  +  +  +  +  +  +  +  +  +
\begin{acknowledgements}
The work of Grimm was supported in part by the NSF through Grant DMS-1217156; 
the work of Beattie was supported in part by the Einstein Foundation - Berlin; 
the work of Gugercin was supported in part by NSF through Grant   DMS-1522616; 
the work of Drma\v{c} was  supported in part by  the Croatian Science Foundation through grant HRZZ-9345.
\end{acknowledgements}

\bibliographystyle{spmpsci}      % mathematics and physical sciences
\bibliography{library}   % name your BibTeX data base

\end{document}